\pdfoutput=1
\documentclass[11pt]{article}
\usepackage[letterpaper,hmargin=2.5cm]{geometry}
\usepackage{amsmath,amssymb,latexsym}
\usepackage{theorem} 
\usepackage{overpic}
\usepackage{esint}
\usepackage{pifont}
\newtheorem{thm}{Theorem}

\newtheorem{lemma}[thm]{Lemma}
\newtheorem{proposition}[thm]{Proposition}
\newtheorem{corollary}[thm]{Corollary}
{\theorembodyfont{\rmfamily} 
\newtheorem{remark}[thm]{Remark}
\newcommand{\sgn}{\mathop{\rm sgn}}

\newcommand{\field}[1]{\mathbb{#1}}
\newcommand{\R}{\field{R}}
\newcommand{\Z}{\field{Z}}
\newcommand{\N}{\field{N}}
\newcommand{\C}{\field{C}}

\newcommand{\CC}{{\mathcal C}}
\newcommand{\OO}{{\mathcal O}}

\renewcommand{\Re}{\mathop{\rm Re}}
\renewcommand{\Im}{\mathop{\rm Im}}

\newcommand{\isdef}{\stackrel{\text{\tiny def}}{=}}

\DeclareRobustCommand{\qed}{%
\ifmmode 
\else \leavevmode\unskip\penalty9999 \hbox{}\nobreak\hfill \fi
\quad\hbox{\qedsymbol}}
\newcommand{\openbox}{\leavevmode
\hbox to.77778em{%
\hfil\vrule
\vbox to.675em{\hrule width.6em\vfil\hrule}%
\vrule\hfil}}
\newcommand{\qedsymbol}{\openbox}
\newcommand{\proofname}{Proof}
\newenvironment{proof}[1][\proofname]{\par
\normalfont \trivlist
\item[\hskip\labelsep   \itshape #1. ]
\ignorespaces
}{%
\qed\endtrivlist }

\begin{document}

\title{Asymptotics of orthogonal polynomials for a weight with a jump on $[-1,1]$}

\author{A.~Foulqui\'{e} Moreno, A.~Mart\'{\i}nez-Finkelshtein, and V.L. Sousa\footnote{Corresponding author.}}
\date{}

\maketitle

\begin{abstract}
We consider the orthogonal polynomials on $[-1,1]$ with respect to the
weight
$$
w_c\left(  x\right)  =h\left(  x\right)  \left(
1-x\right)  ^{\alpha}\left(  1+x\right)  ^{\beta}\Xi_{c}\left(  x\right) , \quad \alpha , \beta >-1,
$$
where $h$ is real analytic and strictly positive on $[-1, 1]$, and $\Xi_{c}$ is a step-like function: $\Xi_{c}(x)=1$ for $x\in [-1, 0)$ and $\Xi_{c}(x)=c^2$,  $c>0$, for $x\in [0, 1]$.
We obtain strong uniform asymptotics of the monic orthogonal polynomials in $\mathbb{C}$, as well as first terms of the asymptotic expansion of the main parameters (leading coefficients of the orthonormal polynomials and the recurrence coefficients) as $n\to \infty$. In particular, we prove for $w_c$ a conjecture of A.~Magnus regarding the asymptotics of the recurrence coefficients. The main focus is on the local analysis at the origin. We study the asymptotics of the Christoffel-Darboux kernel in a neighborhood of the jump and show that the zeros of the orthogonal polynomials no longer exhibit clock behavior.

For the asymptotic analysis we use the steepest descendent method of Deift and Zhou applied to the non-commutative Riemann-Hilbert problems characterizing the orthogonal polynomials. The local analysis at $x=0$ is carried out in terms of confluent hypergeometric functions. Incidentally, we establish some properties of these functions that may have an independent interest.
\end{abstract}




\section{Introduction and statement of results}

\subsection{Introduction} \label{sec:intro}

Szeg\H{o} is the founder of the modern asymptotic theory of orthogonal polynomials on the unit interval for weights $w$ that satisfy the Szeg\H{o} condition
\begin{equation} \label{SzegoCondition}
        \int_{-1}^1 \frac{\log w(x)}{\sqrt{1-x^2}} dx > - \infty.
    \end{equation}
For the classical Jacobi weights the asymptotic results both on and away from the interval of orthogonality, as well as at its endpoints, can be derived using multiple identities that these orthogonal polynomials satisfy: the differential equation, the Rodrigues formula, integral representation, etcetera. However, in a general situation the problem is much more difficult. Starting from the 80's, many new asymptotic results were found for various classes of weights, and the breakthrough was partially motivated by the development of the tools from potential theory and operator theory.

An important new technique  for obtaining asymptotics
for orthogonal polynomials in all regions of the complex plane is based on the characterization of
the orthogonal polynomials by means of a Riemann--Hilbert problem
for $2\times 2$ matrix valued functions due to Fokas, Its, and
Kitaev \cite{Fokas92}, combined with the
steepest descent method of Deift and Zhou, introduced in \cite{MR94d:35143} and further developed in \cite{MR2001m:05258a, MR98b:35155, MR96d:34004}, to mention a few.

A crucial contribution to this method is \cite{MR2087231}, where the complete asymptotic expansion for the orthogonal polynomials with respect to a Jacobi weight modified by a real analytic and strictly positive function is obtained. However, not much is known in the case when the weight has a jump discontinuity on the interval. So far, the only contribution is \cite{Its07b}, where the authors considered an exponential weight on $\R$ with a jump at the origin, although from a different perspective of asymptotics of Hankel determinants.

Combining ideas from \cite{Its07b} and \cite{MR2087231}, we consider polynomials that are orthogonal on a finite interval $[-1,1]$ with respect to a modified Jacobi weight with a jump, namely
    \begin{equation} \label{Definitiew}
        w_c(x) = (1-x)^{\alpha} (1+x)^{\beta} h(x)\, \Xi_{c}(x), \qquad x \in [-1,1],
    \end{equation}
where $\alpha, \beta > -1$ and $h(x)$ is real analytic and
strictly positive on $[-1,1]$, and $ \Xi_{c}$ is a step-like function, equal to $1$ on $[-1,0)$ and $c^2>0$ on $[0,1]$. Observe that $w_1$ is the weight considered in \cite{MR2087231}. The main asymptotic difference between the polynomials orthogonal with respect to $w_1$ and $w_c$, ($c\neq 1$), lies in their behavior near the origin. While in both cases the analysis near the endpoints of the
interval typically involves Bessel functions, only for $c\neq1$ do confluent hypergeometric functions appear around the origin.

We use $P_n(x) = P_n(x;w_c)$ to denote the monic polynomial of degree $n$
orthogonal with respect to the weight $w_{c}$ on $[-1,1]$,
    \[
        \int_{-1}^1 P_n(x;w_c) x^k w_c(x) \ dx = 0, \qquad
        \mbox{for $k = 0,1, \dots, n-1$,}
    \]
and $p_n(x) = p_n(x;w_c)$ to denote the corresponding orthonormal polynomials,
    \[
        p_n(x) = k_n P_n(x),
    \]
where $k_n > 0$ is the leading coefficient of $p_n$.

The leading term of the asymptotics of polynomials $p_n$ and $P_n(z)$ (as $n \to \infty$)  for a weight satisfying the Szeg\H{o} condition \eqref{SzegoCondition} (and $w_c$ does) and $z\in \C\setminus [-1,1]$ is well known, see \cite{szego:1975}. It can be formulated in terms of two functions that will play a relevant role in what follows, and that we introduce here. Namely,
\begin{equation}\label{phi}
    \varphi\left(  z\right)  =z+\sqrt{z^{2}-1}%
\end{equation}
is the conformal map from $\mathbb C \setminus [-1,1]$ onto
the exterior of the unit circle, with the branch of $\sqrt{z^2-1}$ that is analytic in $\mathbb C \setminus [-1,1]$
and behaves like $z$ as $z \to \infty$. Furthermore, since the weight $w_c$ on $[-1,1]$ satisfies \eqref{SzegoCondition}, we can define the so-called Szeg\H{o} function $D(z) = D(z;w_c)$ associated with $w_c$, given by
    \begin{equation*} 
        D(z)=\exp\left( \frac{\sqrt{z^{2}-1}}{2\pi}\int_{-1}^{1}
        \frac{\log w_c(x)}{\sqrt{1-x^{2}}} \frac{dx}{z-x}\right),\qquad\mbox{for
        $z\in\mathbb{C}\setminus[-1,1]$,}
    \end{equation*}
again with $\sqrt{z^2-1}>0$ for $z>1$ and $\sqrt{1-x^{2}}>0$ on $(-1,1)$.
The function $D(z)$ is a non-zero analytic function on $\mathbb C \setminus [-1,1]$
such that
\[ D_+(x) D_-(x) = w_c(x), \qquad \mbox{ for a.e.\ } x \in (-1,1), \]
where $D_+(x)$ and $D_-(x)$ denote the limiting values of $D(z)$ as $z$ approaches $x$ from above and
below, respectively. In particular, by \eqref{SzegoCondition}, the limit
\[ D_{\infty} = \lim_{z \to \infty} D(z) =
    \exp \left( \frac{1}{2\pi} \int_{-1}^1 \frac{\log w_c(x)}{\sqrt{1-x^2}} dx \right) \]
exists and is a positive real number. From Szeg\H{o}'s theory (see \cite{szego:1975}) it follows that
\begin{equation}\label{basicSzego}
    \frac{2^n P_n(z)}{\varphi(z)^{n}} =
        \frac{D_{\infty}}{D(z; w_c)} \frac{\varphi(z)^{1/2}}{\sqrt{2}(z^2-1)^{1/4}}\,
        \left[ 1 +  o(1)\right], \qquad \mbox{as $n\to\infty$,}
\end{equation}
    uniformly on compact subsets of $\overline{\mathbb C} \setminus [-1,1]$. Using the multiplicative property of the Szeg\H{o} function, we conclude that in comparison with the case $c=1$, for $c\neq 1$  there is an extra factor, corresponding to the Szeg\H{o} function of the pure jump $\Xi_{c}$.

In this paper, we give  uniform and more precise asymptotic results for the special weights
(\ref{Definitiew}). We obtain the first terms of the asymptotic
expansions for $k_n$, $P_n$, and $p_n$, as well as for the
coefficients $a_n$ and $b_n$ in the three-term recurrence relation
    \begin{equation} \label{RecRelphin}
        P_{n+1}(z) = (z-b_n) P_n(z) - a_n^2 P_{n-1}(z),
    \end{equation}
satisfied by the monic orthogonal polynomials.

From our analysis we are also able to derive strong asymptotics
for the orthogonal polynomials in the open interval $(-1,1)$, near the endpoints $\pm 1$, and what is most interesting, in a neighborhood of the origin where the jump of the weight takes place.

Since the behavior of the polynomials $P_n$ away from the origin is very similar to the case $c=1$ treated in  \cite{MR2087231}, we will not present all formulas here. However, all the ingredients are contained in the results of the steepest descent analysis performed in Section \ref{sec:RHanalysis}, so that an interested reader can effortlessly derive the omitted asymptotic formulas. In this paper we concentrate on the features of the polynomials and their coefficients that stem from the discontinuity of the weight at the origin.

\subsection{Asymptotics away from the interval of orthogonality}

In order to formulate our results we need to introduce some notation. For $h(x)$ real analytic and
strictly positive on $[-1,1]$, $n\in \N$ and $c>0$ we define the following real-valued function and real quantities:
\begin{align}
\label{def_hbar1}
    \hbar(x) & \isdef \frac{\sqrt{1-x^2}}{2 \pi  }\fint_{-1}^{1}\frac{\log h\left(  t\right)
}{\sqrt{1-t^{2}}}\frac{dt}{t-x}, \qquad x \in (-1,1),\\
\label{defOfEta}
\eta_n & =\eta_n (c)\isdef  \frac{\log c}{\pi}\, \log(4n) + \frac{n \pi}{2} + \frac{\beta -\alpha }{4}\, \pi + \hbar(0),
\end{align}
where $\fint$ is the integral understood in terms of its principal value. In general, we assume always $\sqrt{1-x^2}>0$ for $x\in(-1,1)$, unless stated otherwise.

We also introduce what will play the role of the main phase shift in all asymptotic formulas,
\begin{align}
\label{defOfTheta} \theta_n & = \begin{cases}
\theta_n (c)\isdef  2\left ( \eta_n -    \arg\left(  \Gamma\left(i\frac{\log c}{\pi} \right) \right)\right), & \text{if } c \neq 1, \\
 2\left ( \eta_n    +\frac{\pi}{2}\right), & \text{if } c = 1,
\end{cases}
\end{align}
where $\Gamma(\cdot)$ is the Gamma function; for purely imaginary values of $\lambda\neq 0$, we take $\arg(\Gamma(\lambda ))\in (-\pi/2, \pi/2)$.

The simplest asymptotic result concerns the monic orthogonal polynomials $P_n$.
Observe that a full asymptotic expansion for the usual Jacobi polynomials ($h \equiv 1$, $c=1$) can be found in  \cite[Theorem 8.21.9]{szego:1975}, while for general real analytic and positive $h$ (but with $c=1$) it was established in \cite{MR2087231}. Here we
find only the first two terms of the asymptotic expansion, improving \eqref{basicSzego}:
\begin{thm} \label{thm:outer}
    We have that
    \[
        \frac{2^n P_n(z)}{\varphi(z)^{n}} =
        \frac{D_{\infty}}{D(z;w_c)} \frac{\varphi(z)^{1/2}}{\sqrt{2}(z^2-1)^{1/4}}
        \left[ 1 +  \frac{\mathcal H_n(z)}{n} +\OO\left( \frac{1}{n^2}\right)
        \right], \qquad \mbox{as $n\to\infty$,}
    \]
    uniformly on compact subsets of $\overline{\mathbb C} \setminus [-1,1]$. The function  $\mathcal H(z)$
    is analytic on $\mathbb C \setminus [-1,1]$, and given by
    \begin{equation} \label{ValueofH}
       \mathcal H_n\left(  z\right)  = -\frac{4\alpha^2-1 }{8(\varphi(z)-1) } + \frac{4\beta^2-1}{8(\varphi(z)+1) } - \frac{ \log (c)   }{  2\pi z \varphi(z)}\, \left( \cos(\theta_n) \varphi(z) +\sin(\theta_n)  -\frac{\log (c) }{\pi}\right),
    \end{equation}
with $\theta_n$ defined in \eqref{defOfTheta}.
\end{thm}
\begin{remark}
A more detailed analysis of the Szeg\H{o} function $D(\cdot; w_c)$ is carried out in Section \ref{sec_Szego}.
We can simplify notation in the formula above observing that
$$
 \frac{(z^2-1)^{1/4}}{\varphi(z)^{1/2}}
$$
is the Szeg\H{o} function for the weight $\sqrt{1-x^2}$ on $[-1,1]$, and it takes the value $2^{-1/2}$ at infinity. Hence,
$$
\frac{D_{\infty}}{D(z;w_c)} \frac{\varphi(z)^{1/2}}{\sqrt{2}(z^2-1)^{1/4}}=\frac{D(\infty; \widehat w_c)}{D(z; \widehat w_c)},
$$
where $\widehat w_c(x)=\sqrt{1-x^2} \, w_c(x)$ is known as the trigonometric weight associated to $w_c$.
\end{remark}


The RH analysis performed below for $z\notin [-1,1]$ allows also to establish a result for some relevant parameters associated with the orthogonal polynomials. Recall that the monic polynomials $P_n$ satisfy the three term recurrence
relation \eqref{RecRelphin}. The asymptotic behavior  of these recurrence coefficients (as $n\to \infty$) is given in the following theorem:
\begin{thm} \label{theoremRecCoef}
    As $n\to \infty$,
  \begin{equation}\label{asymptotics A}
    a_{n}      =\frac{1}{2}-\frac{\log c}{ 2 \pi n}\sin (\theta_{n})+\OO\left(  \frac{1}{n^{2}%
}\right),
\end{equation}
    \begin{equation}
    b_n = -\frac{\log c}{  \pi n}\,   \cos \left(\theta_n \right)   +\OO\left(  \frac{1}{n^{2}}\right), \label{asymptotics B}
    \end{equation}
with $\theta_n$ defined in \eqref{defOfTheta}.
\end{thm}
\begin{remark}
In  \cite{Magnus1995}, A.~Magnus studied weights of the form
\begin{equation*}
  \left(  1-x\right)  ^{\alpha} \left(  1+x\right)  ^{\beta}\left|
x_{0}-x\right|  ^{\gamma} \times
\begin{cases}
B, & \text{for }x\in\left[  -1,x_{0}\right)
\text{,}\\
A, & \text{for }x\in\left[  x_{0},1\right]  \text{,}%
\end{cases}
\end{equation*}
with $A$ and $B>0$ and $\alpha$, $\beta$ and $\gamma>-1$, and $x_{0}\in\left(
-1,1\right)  $. Formulas \eqref{asymptotics A}--\eqref{asymptotics B} show that for $\gamma=x_0=0$ the asymptotic behavior of the recurrence coefficients conjectured in \cite{Magnus1995}  is correct, with the possibility to replace $o(1/n)$ by $\OO(1/n^2)$ in the error term. For more details see Section \ref{sec:asymptRecurrencproofs} below; the proof of the conjecture in its full generality is contained in \cite{FMS2}.
\end{remark}

The leading coefficients $k_n$ of the orthonormal polynomials $p_n$ satisfy the following asymptotic relation:
\begin{thm}
\label{thm:leading} As $n\rightarrow\infty$,
\[
k_{n}=\frac{2^{n}}{\sqrt{\pi}D_{\infty}}\left[  1- \left(\frac{2 \alpha^2+ 2\beta^2 -1 }{8  }    + \frac{\log(c)}{2 \pi }\, \left( \frac{\log(c)}{ \pi } + \sin(\theta_{n+1})  \right)  \right)\,  \frac{1}{  n }   + \OO\left( \frac{1}{n^2}\right)  \right]  \text{,}%
\]
with $\theta_n$ defined in \eqref{defOfTheta}.
\end{thm}

\subsection{Local asymptotics}

Now we need to introduce further notation. Set
\begin{equation}\label{defGIntro}
    G\left(  a;\zeta\right)  \isdef  {_{1}F_{1}}\left(  a;1;\zeta\right) e^{-\zeta/2}=e^{-\zeta/2}\, \sum_{k=0}^\infty \frac{(a)_k}{(k!)^2}\, \zeta^k,
\end{equation}
where ${_{1}F_{1}}\left(  a;b;\cdot\right)$ is the confluent hypergeometric function; $G$ is an entire function of $\zeta$ for any value of the parameter $a\in \C$, and $G(a; 0)=1$.
Furthermore, for $x\in (-\delta, 0)\cup(0, \delta)$ let
\begin{equation}\label{defFuncRho}
  \rho (x)\isdef   \dfrac{\log c}{\pi} \,  \log\left|   \dfrac{\arcsin(x) }{2 x}\, \left(1+ \sqrt{1-x^{2}} \right)\right| -\frac{\alpha +\beta }{2}\,\arcsin(x)+ \hbar(x)-\hbar(0),
\end{equation}
completed to a continuous function on $(-1,1)$ by $\rho(0)=0$.

Set also
\begin{equation}\label{defUpsilon}
\Upsilon(c)\isdef \sgn(\log(c))\, \sqrt{\frac{2 c \log c}{c^{2}-1}}, \quad c\neq1, \quad \Upsilon(1)=1,
\end{equation}
where we always take the positive value of the square root.
The asymptotic behavior of $P_n$ on compact subsets of an interval $(-\delta , \delta )\subset (-1,1)$ is given by the following theorem:
\begin{thm}
\label{thm:local} For $\delta\in (0,1)$, locally uniformly on compact subsets of $(-\delta , \delta )$ the following asymptotic formula holds:
\begin{equation*}
\begin{split}
  P_n(x)    &=\frac{D_\infty}{2^{n-1/2}\sqrt{c\,  w_1(x)}} \frac{\Upsilon( c)}{  (1-x^2)^{1/4}} \\ &\times \Re\left[e^{i(\rho (x)+ \frac {\theta_n -\pi-\arcsin(x)}{2})} G\left(  \lambda ;2 i n \arcsin\left(  x\right)\right) \left( 1 +  \frac{ \mathcal R_n(x)}{n} +\OO\left( \frac{1}{n^2}\right)\right)  \right],
\end{split}
\end{equation*}
with
\begin{equation*}%
\begin{split}
 \mathcal{R} _{n}\left(  x\right)    & = -\frac{4\alpha^2-1 }{8( e^{  i\arccos (x) }-1) } + \frac{4\beta^2-1}{8( e^{  i\arccos (x) }+1) } \\ & - \frac{ \log (c)   }{  2\pi x  e^{  i\arccos (x) }}\,   \left( \cos(\theta_n)e^{  i\arccos (x) } +\sin(\theta_n)  -\frac{\log (c) }{\pi}\right)  \\ & +
\frac{i \log c}{2\pi \arcsin(x)  }\left( \frac{\log c}{\pi }+ e^{-i\,(2\rho(x)+\theta
_{n}+ \arccos(x))}  \right)  ,
\end{split}
\end{equation*}%
$\rho(x)$ given in \eqref{defFuncRho}, $\theta_n$ in \eqref{defOfTheta}, and $\lambda =i\log(c)/\pi$.
\end{thm}
\begin{remark}
The Riemann-Hilbert analysis we perform next gives us an asymptotic expression for $P_{n}$'s in a small disk of the complex plane centered at the origin, see formula \eqref{expressioncloseto0} in Section \ref{sec:asymptoninterval}.
\end{remark}
\begin{corollary}\label{cor:localbeh}
Locally uniformly for $x\in (-\delta , \delta )$, $\delta\in (0,1)$,
\begin{equation}\label{asymptPicloseto0scaled}
\begin{split}
  P_n\left(\frac{\pi x}{n}\right) & =\frac{D_\infty\, \Upsilon(c)}{2^{n-1/2}\sqrt{c\, h(0)}}\,  \Im \left[   e^{ i \,  \theta_n/2} \,     G\left(  \lambda ;2 \pi i x\right) \, \left( 1 +   \OO\left( \frac{1}{n }\right)\right)  \right],
\end{split}
\end{equation}
with  $\theta_n$ given in \eqref{defOfTheta} and $\lambda =i\log(c)/\pi$.
\end{corollary}

\begin{figure}[htb]
\centering \begin{overpic}[scale=1]%
{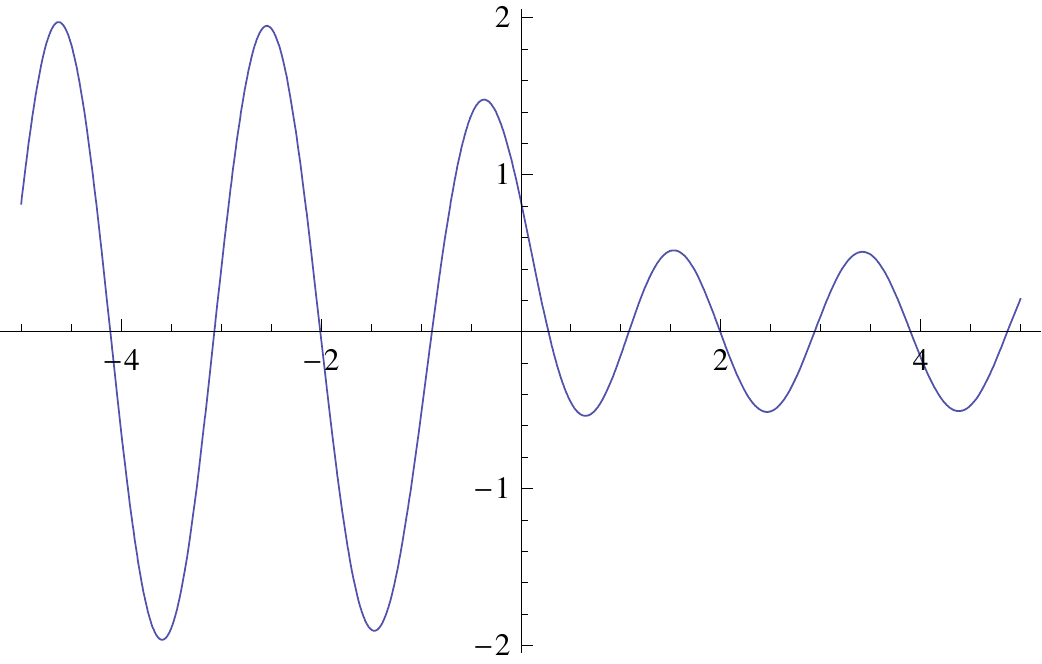}%
\end{overpic}
\caption{Typical graphics of the r.h.s.\ of \eqref{asymptPicloseto0scaled} near the origin.}
\label{fig:localzeros}
\end{figure}

See Figure \ref{fig:localzeros} for a typical behavior of the function in the right hand side of \eqref{asymptPicloseto0scaled} close to the origin.

Recall that $P_{n}$ has $n$ simple zeros, all lying on $(-1,1)$. It is well known that they distribute asymptotically in the weak-* sense according to the equilibrium measure of the interval. In other words, the normalized zero counting measure for the sequence $P_{n}$ weakly tends to the absolutely continuous measure on $[-1,1]$ given by $\omega(x)\, dx$, with
$$
 \omega(x)\isdef \frac{1}{\pi}\, \frac{1}{\sqrt{1-x^{2}}}.
$$
As it follows from several works of Deift and collaborators (and also from a recent series of papers of Lubinsky and Levin and Lubinsky, see e.g.~\cite{Lubinsky2009,Lubinsky2008}), a much stronger statement holds: at any point of $(-1,0)\cup(0,1)$ they distribute very precisely in accordance with $\omega(x)$, complying with the so-called ``clock behavior'', see e.g.~\cite{Simon:2009ly}. If, following \cite{Simon:2009ly}, we enumerate the zeros $x_{j}^{(n)}$ of $P_{n}$ as follows,
\begin{equation}
\label{enumerationzerosX}
\dots <x_{-k}^{(n)} <\dots < \dots <x_{-1}^{(n)} <0\leq x_{0}^{(n)}<\dots<x_{k}^{(n)}<\dots
\end{equation}
then ``clock behavior'' at the origin (where $\omega(0)=1/\pi$) means
\begin{equation}\label{clockbehavior}
\lim_{n\to \infty} \frac{n}{\pi}\, \left(x_{j+1}^{(n)}-x_{j}^{(n)} \right)=1, \quad j\in \Z.
\end{equation}

\begin{proposition}
\label{prop:clockbehavior}
If $c > 1$, then the sequence $\{n\, x_{0}^{(n)}/\pi\}$ is dense in an interval of the form $[0, t]$, where $t=t(c)<1$. Furthermore,
$$
0<\liminf_{n} \frac{n }{\pi} \left( x^{(n)}_{k} -  x_{k-1}^{(n)} \right)
\leq \limsup_{n}\frac{n }{\pi} \left(x^{(n)}_{k} -  x_{k-1}^{(n)} \right)<1, \quad k \in \N,
$$
and
$$
\liminf_{n} \frac{n }{\pi} \left( x^{(n)}_{k} -  x_{k-1}^{(n)} \right)>1, \quad -k \in \N.
$$
In particular, the clock behavior of the zeros of $P_{n}$ at the origin does not hold.

If $0<c<1$, the same inequalities hold inverting the roles of $k$ and $-k$.
\end{proposition}
This result is not surprising, taking into account that $x=0$ is not even a Lebesgue point for the weight $w_{c}$, that is,  regardless of the meaning we give to $w_{c}(0)$,
$$
\lim_{s \to 0^{+}} \frac{1}{s}\, \int_{-s}^{ s} \left| w_{c}(x)-w_{c}(0) \right|\, dx \neq 0.
$$
However, to the best of our knowledge, $w_{c}$ with $c\neq 1$ provides the first instance of an explicit orthogonality measure for which the clock behavior fails in the bulk (interior of its support).

\begin{remark} \label{rem:quasi}
A weaker condition than \eqref{clockbehavior} is the quasi-clock behavior (see \cite{Simon:2009ly}), namely
$$
\lim_{n\to \infty} \frac{ x_{j+1}^{(n)}-x_{j}^{(n)} }{x_{1}^{(n)}-x_{0}^{(n)} }=1, \quad j\in \Z.
$$
This limit is violated in our situation too. However,
\begin{equation}\label{clockinthelimit}
\lim_{j \to \pm\infty} \liminf_{n\to \infty} \frac{n}{\pi}\, \left(x_{j+1}^{(n)}-x_{j}^{(n)} \right) =\lim_{j \to \pm\infty} \limsup_{n\to \infty} \frac{n}{\pi}\, \left( x_{j+1}^{(n)}-x_{j}^{(n)} \right) =1,
\end{equation}
which shows a smooth transition to the genuine clock behavior as we move away from the jump of the weight.
\end{remark}

Very much related with the clock behavior is the ``universality problem'' for the Christoffel-Darboux (or CD) kernel
\begin{equation}\label{defCDKernel}
    K_{n}\left(  x,y\right)  \isdef \sum_{k=0}^{n-1}p_{k}\left(  x\right)  p_{k}\left(
y\right)  ,
\end{equation}
where $  p_{n}  $ are the orthonormal
polynomials with respect to the weight $w_c  $. This problem has its origin  in the random matrix theory and has been attracting lately close attention of many researchers. A recent series of remarkable contributions of Lubinsky allowed to weaken considerably the conditions on the weight to be able to assure universality: now we know that for $t$ within the support of the weight where it is continuous,
\begin{equation}\label{kernelsine}
\lim_{n\to \infty}\frac{\pi}{n \sqrt{1-t^2}}\, K_{n}\left(t+  \frac{\pi x }{n \sqrt{1-t^2}},t+ \frac{\pi y }{n\sqrt{1-t^2}}\right) =\frac{\sin \left(\pi(x-y)\right)}{\pi( x-y) } .
\end{equation}
The right hand side is the well-known sine (or sinc) kernel; for our weight $w_c$, this formula is valid for $t\in (-1,0)\cup(0,1)$. It was observed in \cite{Lubinsky2008} that \eqref{kernelsine} implies \eqref{clockbehavior}.

We show that the jump discontinuity in the weight leads to a different kernel, constructed in terms of the confluent hypergeometric function defined in \eqref{defGIntro}:
\begin{thm}
\label{thm:kernel} For $c>0$, $c\neq 1$, locally uniformly for $x$ and $y$ on $(- \delta, \delta)$, $0<\delta<1$, 
\begin{equation} \label{kernelFinal1}
\lim_{n\to \infty}\frac{\pi}{n}\, K_{n}\left(  \frac{\pi x }{n},\frac{\pi y }{n}\right)  = K_{\infty}\left( x,y \right)  ,
\end{equation}
with
\begin{equation} \label{limitkernel}
\begin{split}
& K_{\infty}\left( x,y \right)    \\ & = \begin{cases}
\dfrac{1}{h(0)\,  \pi i} \dfrac{  \log c}{c^2-1} \, \dfrac{\big[  G\left(  1+\lambda ;2 \pi i x \right); G\left(   \lambda ;2 \pi i y \right) \big]}{ x-y }, & x\neq y, \\[3mm]
\dfrac{2}{h(0)}  \dfrac{  \log c}{c^2-1}    \left(  G'\left(  1+\lambda ;2 \pi i x \right) G\left(   \lambda ;2 \pi i x \right)-  G\left(  1+\lambda ;2 \pi i x \right) G'\left(   \lambda ;2 \pi i x \right) \right), & x=y,
\end{cases}
\end{split}
\end{equation}
where $\lambda =i\log(c)/\pi$, $G$ was introduced in \eqref{defGIntro}, and as usual, $[f(x); g(y)]=f(x) g(y)-f(y)g(x)$.
\end{thm}

Several remarks are in order.

Since $G'\left(  1+\lambda ;0 \right) =\lambda +1/2$ and $G\left(   \lambda ;0 \right)= 1  $, evaluating $K_{\infty}(0,0)$ in \eqref{limitkernel} we conclude that
$$
\lim_{n\to \infty} \frac{K_{n}\left(  \frac{\pi x }{n},\frac{\pi x }{n}\right) }{K_{n}\left( 0,0\right) }=G'\left(  1+\lambda ;2 \pi i x \right) G\left(   \lambda ;2 \pi i x \right)-  G\left(  1+\lambda ;2 \pi i x \right) G'\left(   \lambda ;2 \pi i x \right),
$$
locally uniformly in $(-\delta,\delta)$. This shows that even the weak Lubinsky's ``wiggle condition'' (term coined by B.\ Simon, see e.g.~\cite[Theorem 3.6]{Simon:2009ly})  is not satisfied in a neighborhood of the jump of the weight.

The kernel for $x\neq y$ in \eqref{limitkernel} is written in the so-called integrable form. Taking into account the properties of the functions in the right hand side, we can rewrite it alternatively in a totally real form:
\begin{equation} \label{kernelFinal2}
\begin{split}
K_{\infty}(x,y) & =  \frac{2}{\pi( x-y) \, h(0)} \frac{  \log c}{c^2-1}  \Im \left( G\left(  1+\lambda ;2 \pi i x \right) G\left(   \lambda ;2 \pi i y \right)   \right), \quad x\neq y.
\end{split}
\end{equation}

Since $G(1,z)=\exp(z/2)$, straightforward computations show that as $c\to 1$, $K_{\infty}$ reduces to the sine kernel.
Notice that combining  ideas from \cite{Lubinsky2008} and \cite{Simon2008} we can use  \eqref{kernelFinal2} to arrive at the same conclusions about the spacing of zeros of $P_n$'s as we did at the end of Subsection \ref{sec:asymptoninterval}.

The confluent hypergeometric functions appeared in the scaling limit (as the number of particles goes to infinity) of the correlation functions of the pseudo-Jacobi ensemble in  \cite{Borodin:2001xr}. This ensemble corresponds to a sequence of weights of the form
\begin{equation}\label{borodin}
(1+x^{2})^{-n-\Re(s)} e^{2\Im(s) \arg(1+i x)}, \quad x\in \R,
\end{equation}
where $n$ is the degree of the polynomial and $s$ is a complex parameter. The connection between both problems becomes apparent if we perform the inversion $x\mapsto 1/x$ in \eqref{borodin}; this creates at the origin an algebraic singularity with the exponent $\Re(s)$ and a jump depending on $\Im(s)$. $K_{\infty}$ is a particular case of the reproducing kernel obtained by Borodin and Olshansky in Theorem 2.1 of \cite{Borodin:2001xr} when $\Re(s)=0$; for a general situation, see  \cite{FMS2}.

A recent paper of Lubinsky \cite{Lubinsky2009} revealed an interesting connection of $K_{\infty}$ with the theory of entire functions. Namely, in accordance with Theorem 1.6 of \cite{Lubinsky2009}, $K_{\infty}$ is a reproducing kernel of a de Brange space, equivalent to a classical Paley-Wiener space.
More precisely and following the notation of \cite{Lubinsky2009}, the \emph{Hermite-Biehler class}  $\overline{HB}$  is the set of entire functions $E$ with no zeros in the upper half plane $\C^{+}\isdef \{\Im z >0\}$ and such that $|E(z)|\geq |E(\overline{z})|$ for $z\in \C^{+}$. The \emph{de Branges space} $\mathcal H(E)$ corresponding to $E\in \overline{HB}$ is comprised of entire functions $g$ such that both $g(z)/E(z)$ and $\overline{g(\overline{z})}/E(z)$ belong to the Hardy class $H^{2}(\C^{+})$. A reproducing kernel for $\mathcal H(E)$ is
\begin{equation}
\label{debrangesReproducing}
\mathcal K(x,y)=\frac{i}{2\pi}\, \frac{E(x) \overline{E(y)}- \overline{E(\overline{x})} E(\overline{y})}{x-\overline{y}}, \quad x\neq y.
\end{equation}
Comparing this expression with $K_{\infty}$ in \eqref{limitkernel} we conclude that 
$$
\mathcal K(x,\overline{y})=K_{\infty}(x,y),
$$
with $\lambda =i\log(c)/\pi$ and
$$
E(z) = \left( \frac{2}{h(0)}\, \frac{\log c}{c^{2}-1} \right)^{1/2} G(\lambda, 2 \pi i z) \in \overline{HB}
$$
(see below). Lubinsky showed that reproducing kernels, different from the right hand side in \eqref{kernelsine}, can appear for sequences of measures (cf.~\cite{Borodin:2001xr}). To the best of our knowledge, this is the first explicit example of a non-sine reproducing kernel of a de Brange space that arises as a universality limit in the bulk of a fixed measure of orthogonality.

The proof of the asymptotic results stated in this paper (see Section~\ref{sec:proofs}) is based on the steepest descent analysis of the Riemann-Hilbert problem that we carry out in Section~\ref{sec:RHanalysis}. A key step is the construction of the local representation at the origin, which is done in Subsection~\ref{subsec:local0}. The study of the zeros of $P_{n}$'s at the origin, the analysis of the clock behavior and the connection with the de Brange spaces requires some further properties of the confluent hypergeometric function ${_{1}F_{1}}\left(  \lambda;1;z\right)$, which we were unable to find in the literature and which might have an independent interest. We summarize them in the next proposition; the proofs are relegated to Section \ref{sec:confluentproperties}. 

\begin{proposition}\label{prop:nozeros}
Let $a\in \R\setminus \{0\}$. Then
\begin{enumerate}
\item[(i)] functions
$$
f_{1}(z)=G(i a,   i z) \quad \text{and} \quad f_{2}(z)=\overline{G(1+i a,   i \overline{z})}
$$
(see \eqref{defGIntro}) belong to the Hermite-Biehler class $\overline{HB}$;
\item[(ii)] for $x \in \R$,
\begin{equation}
\label{Fnonzero}
 {_{1}F_{1}}\left(  i a ;1;  i x\right) \neq 0 \quad \text{and} \quad  {_{1}F_{1}}\left( 1+ i a ;1;  i x\right) \neq 0.
\end{equation}
In particular, all zeros of $ {_{1}F_{1}}\left(  i a ;1;  i z\right) $ lie in the lower half plane $\C^{-}\isdef \{\Im z < 0\}$, while the zeros of $ {_{1}F_{1}}\left( 1+ i a ;1;  i z\right) $ lie in the upper half plane $\C^{+}$. Additionally,
\begin{equation}
\label{inequalitystrictforF}
|{_{1}F_{1}}\left( 1+ i a ;1;  i z\right)| \leq |{_{1}F_{1}}\left(   i a ;1;  i z\right)|,  \quad \Im z \geq 0,
\end{equation}
and the equality holds only for $z\in \R$.
\item[(iii)] if $a>0$, the function
$$
y(x)\isdef \arg  {_{1}F_{1}}\left(  ia ;1;  i x\right), \quad y(0)=0,
$$
is real-analytic and non-positive, strictly increasing on the negative and strictly decreasing on the positive semiaxis. For $a<0$ the same assertion is valid replacing $y(x)$ by $- y(x)$.
It is also the solution of the following initial value problem:
\begin{equation}\label{3}
x y' =a \left(\cos\left(x-2y\right)-1 \right), \quad y(0)=0;
\end{equation}
\item[(iv)] for $a\in \R$, the function
\begin{equation}
\label{defFrakGbis}
\mathfrak G(x) \isdef x- 2\, \arg \left(  {_{1}F_{1}}\left( ia;1; i x\right)\right)=x-2y(x),\quad \mathfrak G(0)=0,
\end{equation}
is strictly increasing in $\R$.
\end{enumerate}
\end{proposition}
\begin{remark}
The assertion in \emph{(i)} does not imply that $ {_{1}F_{1}}\left(  i a ;1;  i z\right)\in \overline{HB}$, and in general, this is not true.

Interestingly enough, the proof of \emph{(i)} is based on some properties of the Christoffel-Darboux kernel observed by Lubinsky in \cite{Lubinsky2009}. In this sense, the theory of the confluent hypergeometric functions has benefited from the properties of the reproducing kernels. In the opposite direction, the strict inequality in \eqref{inequalitystrictforF} implies that $K_{\infty}(z, \overline{z})>0$ for $z\in \C\setminus \R$ and $c\neq 1$, see \eqref{inequalityForK} below.
\end{remark}

\section{The steepest descent analysis} \label{sec:RHanalysis}

\subsection{The Riemann-Hilbert problem}

Following Fokas, Its and Kitaev \cite{Fokas92} we characterize both the orthogonal polynomials and the CD kernel in terms of the unique solution $\mathbf Y$ of the following $2\times2$ matrix valued
Riemann-Hilbert (RH) problem: for $n\in \N$,
\begin{enumerate}
\item[(Y1)] $\mathbf Y$ is analytic in $\C\setminus [-1,1]$.
\item[(Y2)] On $(-1,0) \cup (0,1)$, $\mathbf Y$ possesses continuous boundary values $\mathbf Y_+$ (from the upper half
plane) and $\mathbf Y_-$ (from the lower half plane), and
\[
\mathbf Y_{+}(x)=\mathbf Y_{-}(x)\,
\begin{pmatrix}
1 & w _c  (x) \\
0 & 1
\end{pmatrix}.
\]
\item[(Y3)] As $z\rightarrow\infty$,
$$
\mathbf  Y(z)=\left(  \mathbf{I}+\OO\left(  \frac{1}{z}\right)  \right)
\begin{pmatrix}
z^{n} & 0\\
0 & z^{-n}%
\end{pmatrix}  ,
$$
where $\mathbf I$ is the identity $2 \times 2 $ matrix.

\item[(Y4)]  $\mathbf Y$ has the following asymptotic behavior at the end points of the interval: for $\zeta\in \{-1, 1\}$ set $s=\alpha$ if $\zeta=1$, and $s=\beta$ if $\zeta=-1$. Then for $z\rightarrow \zeta $, $z\in\mathbb{C}\backslash\left[  -1,1\right]  $,
\[
\mathbf Y(z)= \begin{cases}
\mathcal O
\begin{pmatrix}
1 & \left\vert z-\zeta\right\vert ^{s}\\
1 & \left\vert z-\zeta\right\vert ^{s}%
\end{pmatrix}, & \text{if }s<0;\\
\mathcal O \begin{pmatrix}
1 & \log\left\vert z-\zeta\right\vert \\
1 & \log\left\vert z-\zeta\right\vert
\end{pmatrix},
  & \text{if }s=0;\\
\mathcal O
\begin{pmatrix}
1 & 1\\
1 & 1
\end{pmatrix},
& \text{if }s>0.%
\end{cases}
\]
Furthermore, at the origin $\mathbf Y$ has the following behavior: for $z\rightarrow0$,   $z\in\mathbb{C}\backslash\left[  -1,1\right]  $,
$$
\mathbf  Y ( z) =\mathcal O
\begin{pmatrix}
1 & \log\left\vert z\right\vert \\
1 & \log\left\vert z\right\vert
\end{pmatrix}.
$$
\end{enumerate}
Standard arguments (see e.g.~\cite{MR2087231}) show that this RH problem has a unique solution given by
\begin{equation}
\mathbf  Y\left(  z, n\right)  =\left(
\begin{array}
[c]{cc}%
P_{n}\left(  z\right)  & \CC\left(  P_{n}w_c\right)  \left(  z\right) \\
-2\pi ik_{n-1}^{2}P_{n-1}\left(  z\right)  & -2\pi ik_{n-1}^{2}\CC\left(
P_{n-1}w_c\right)  \left(  z\right)
\end{array}
\right)  \text{,} \label{sol Y}%
\end{equation}
where $P_{n}$ is the monic orthogonal polynomial of degree $n$\ with respect
to $w_c$; $k_{n}$  is the leading coefficient of the orthonormal polynomial $p_{n}$,
and $\CC\left(  \cdot \right)  $\ is the Cauchy transform on $\left[  -1,1\right]  $
defined by
\[
\CC\left(  f\right)  \left(  z\right)  =\frac{1}{2\pi i}\int_{-1}^{1}%
\frac{f\left(  x\right)  }{x-z}\, dx\,.
\]
Clearly, $\mathbf Y$ and other matrices introduced hereafter depend on $n$, fact that we indicate writing $\mathbf  Y(\cdot, n)$. However, we omit the explicit reference to $n$ from the notation whenever it cannot lead us into confusion.

\subsection{First transformations} \label{sec:First transformations}

We apply the Deift-Zhou method of steepest descent to the RH problem above; some of the steps are standard and we occasionally omit those less relevant details that can be easily found in literature (each time we try to provide a suitable reference though). As in \eqref{phi}, $\varphi$ denotes the conformal mapping from $\C \setminus \left[  -1,1\right]  $ onto the exterior of the unit circle.
Let $\sigma_3=\begin{pmatrix}
1 & 0 \\ 0 & -1
\end{pmatrix}$ be the third Pauli matrix; in what follows, for $a\in \C \setminus \{0\}$  we use the notation
$$
a^{ \sigma_{3}} \isdef \begin{pmatrix}
a  & 0\\
0 & 1/a
\end{pmatrix};
$$
then for $b\in \C$, $a^{ b\sigma_{3}}$ is understood as $(a^b)^{ \sigma_{3}}$. Furthermore, if $\gamma$ is an oriented Jordan arc, and an analytic function $f$ has boundary values at $\gamma $, we denote by $f_+$ (resp., $f_-$) its boundary values on $\gamma $ from the left (resp., from the right).

Set
\begin{equation}
\mathbf T\left(  z \right)  \isdef 2^{n\sigma_{3}} \mathbf Y\left(  z \right)  \varphi\left(
z\right)  ^{-n\sigma_{3}}. \label{sol T}%
\end{equation}
Then $\mathbf T$ is the unique solution of the following equivalent RH problem:
\begin{enumerate}
 \item[(T1)] $\mathbf T$ is analytic in $\C\setminus [-1,1]$.
\item[(T2)] On $(-1,0) \cup (0,1)$, oriented from $-1$ to $1$, $\mathbf  T$ possesses continuous boundary values $\mathbf T_+$  and $\mathbf T_-$, and
\[
\mathbf T_{+}(x)= \mathbf T_{-}(x)\,
\begin{pmatrix}
\varphi_{+}^{-2n}\left(  x\right)  & w_c \left(  x\right) \\
0 & \varphi_{-}^{2n}\left(  x\right)
\end{pmatrix}.
\]
\item[(T3)] As $z\rightarrow\infty$,
$$
\mathbf T(z)=  \mathbf{I}+\OO\left(  \frac{1}{z}   \right).
$$

\item[(T4)]  $\mathbf T$ has the same asymptotic behavior as $\mathbf Y$ at $\pm 1$ and $0$ .
\end{enumerate}

Next transformation is based upon the factorization of the jump matrix for $\mathbf T$:
\begin{equation}
\begin{pmatrix}
\varphi_{+}^{-2n} & w_c\\
0 & \varphi_{-}^{-2n}%
\end{pmatrix}
 =\begin{pmatrix}
1 & 0\\
w_c^{-1}\varphi_{-}^{-2n} & 1
\end{pmatrix} \,
\begin{pmatrix}
0 & w_c\\
-1/w_c & 0
\end{pmatrix}\,
\begin{pmatrix}
1 & 0\\
w_c^{-1} \varphi_{+}^{-2n} & 1
\end{pmatrix}. \label{Dec Vt}%
\end{equation}
In order to introduce a contour deformation we need to extend the definition of the weight of orthogonality to a neighborhood of the interval $[-1,1]$.

By assumptions, $h$ is a holomorphic function in a neighborhood $U$ of $[-1,1]$, and positive on this interval. We set
\begin{equation}\label{defwanalytic}
    w(z)\isdef h(z)\, \left(  1-z\right)  ^{\alpha}\left(  1+z\right)  ^{\beta}
\end{equation}
holomorphic in $U\setminus \left( (-\infty, -1] \cup [1, +\infty) \right)$, and such that $w(x)>0$ for $x\in (-1,1)$. In particular, $w(0)=h(0)$. We also extend the definition of the step function $\Xi_{c}$ by
\[
\Xi_{c}\left(  z\right)  =
\begin{cases}
1 , & \text{if } \Re z<0\\
c^2, & \text{if } \Re z\geq 0.
\end{cases}
\]
Then we set
\begin{equation} \label{defwcanalytic}
w_c(z) \isdef w(z)\, \Xi_{c}(z),
\end{equation}
which is a holomorphic function in $\widetilde U\isdef U\setminus \left( (-\infty, -1] \cup [1, +\infty) \cup i\R \right)$.

With this definition the left and rightmost matrices in \eqref{Dec Vt} have an analytic extension to the portion of $\widetilde U$ in the lower and upper half plane, respectively, and we can define the next step: lens opening or contour deformation.
\begin{figure}[htb]
\centering \begin{overpic}[scale=1.5]%
{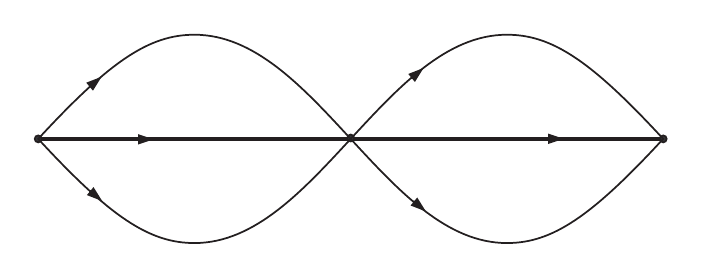}%
     \put(0,20){$\small -1 $}
\put(96,20){$\small 1$}
  \put(49,22){$\small 0 $}
  \put(19,35){$\small \gamma_1$}
    \put(19,3){$\small \gamma_2$}
      \put(78,35){$\small \gamma_3$}
    \put(78,3){$\small \gamma_4$}
    \put(59,14){\small inner lower dom.}
    \put(59,24){\small inner upper dom.}
    \put(41,33){\small outer dom.}
\end{overpic}
\caption{First lens opening.}
\label{fig:lenses1}
\end{figure}
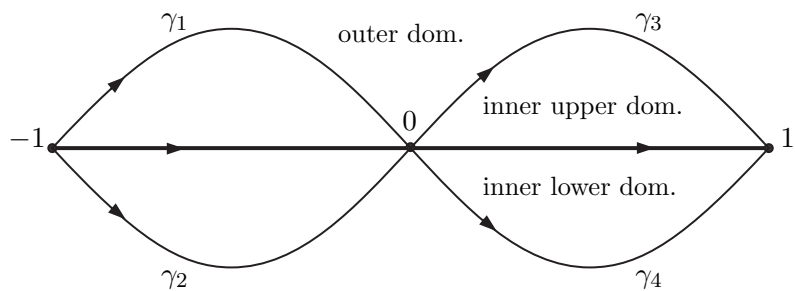
Namely, we build four new contours $\gamma_i$ lying in $\widetilde U$ (except for their end points) such that $\gamma_1$ and $\gamma_3$ are in the upper half plane, and $\gamma_1$ and $\gamma_2$ are in the left half plane, and oriented ``from $-1$ to $1$'' (see Fig.~\ref{fig:lenses1}). This construction defines three domains: the inner upper domain, bounded by $[-1,1]$ and the curves $\gamma_1$ and $\gamma_3$; the inner lower domain, bounded by $[-1,1]$ and the curves $\gamma_2$ and $\gamma_4$, and finally the outer domain, bounded by curves $\gamma_i$ and containing the infinity.

Using the matrix $\mathbf T$ from \eqref{sol T} we define the new matrix  $\mathbf S$  by
\begin{equation}
\mathbf S\left(  z\right)  \isdef
\begin{cases}
\mathbf T\left(  z\right),  & \text{for $z$ in the outer domain, }\\
\mathbf T(z)\, \begin{pmatrix}
1 & 0\\
-\frac{1}{w_c(z)}\varphi^{-2n}(z) & 1
\end{pmatrix},
  & \text{for $z$ in the inner upper domain, }\\[3mm]
\mathbf T(z)\, \begin{pmatrix}
1 & 0\\
\frac{1}{w_c(z)}\varphi^{-2n}(z) & 1
\end{pmatrix},
   & \text{for $z$ in the inner lower domain.}%
\end{cases}
  \label{sol S}%
\end{equation}
Then $\mathbf S$ is the unique solution of the following RH problem:
\begin{enumerate}
\item[(S1)] $\mathbf S$ is analytic in $\mathbb{C}\backslash\Sigma$, where $\Sigma \isdef [-1,1]\cup \bigcup_{i=1}^4\gamma_i$.

\item[(S2)] $\mathbf S$ satisfies the following jump relations:
\begin{align*}
\mathbf S_{+}(z)  &  = \mathbf S_{-}(z)
\begin{pmatrix}
1 & 0\\
\frac{1}{w_c\left(  z\right)  }\, \varphi\left(  z\right)  _{{}}^{-2n} & 1
\end{pmatrix},
   \text{ for } z\in\left( \bigcup_{i=1}^4\gamma_i\right) \setminus \{-1, 0, 1 \},\\
\mathbf S_{+}(x)  &  = \mathbf S_{-}(x)
\begin{pmatrix}
0 & w_c \left(  x\right) \\
-\frac{1}{w_c\left(  x\right)  } & 0
\end{pmatrix}
  \text{, \ \ \ for }x\in(-1,0)\cup (0,1).%
\end{align*}

\item[(S3)] As $z\rightarrow\infty$,
\[
\mathbf S(z)= \mathbf{I}+\OO\left(  \frac{1}{z}\right)  .
\]

\item[(S4)] $\mathbf S$ has the following asymptotic behavior at the end points of the interval: for $\zeta\in \{-1, 1\}$  set $s=\alpha$ if $\zeta=1$, and $s=\beta$ if $\zeta=-1$. Then for $z\rightarrow \zeta  $, $z\in\mathbb{C}\setminus \Sigma  $,
\begin{itemize}
\item for $s<0$:%
\[
\mathbf S\left(  z\right)  =\mathcal O
\begin{pmatrix}
1 & \left\vert z-\zeta\right\vert ^{s}\\
1 & \left\vert z-\zeta\right\vert ^{s}%
\end{pmatrix}
  \text{, \ \ \ as }z\rightarrow\zeta ;%
\]

\item for $s=0$:%
\[
\mathbf S\left(  z\right)  =\mathcal O\left(
\begin{array}
[c]{cc}%
\log\left\vert z-\zeta\right\vert  & \log\left\vert z-\zeta\right\vert \\
\log\left\vert z-\zeta\right\vert  & \log\left\vert z-\zeta\right\vert
\end{array}
\right)  \text{, \ \ \ as }z\rightarrow \zeta;%
\]

\item for $s>0$:%
\[
\mathbf S(z)=
\begin{cases} \mathcal O
\begin{pmatrix}
1 & 1\\
1 & 1
\end{pmatrix},
 & \text{as }z\rightarrow \zeta\text{ from the outer domain;}\\
\mathcal O
\begin{pmatrix}
\left\vert z-\zeta\right\vert ^{-s} & 1\\
\left\vert z-\zeta\right\vert ^{-s} & 1
\end{pmatrix},
  & \text{as }z\rightarrow \zeta\text{ from the inner domains.}%
\end{cases}
\]

\end{itemize}

\item[(S5)] $\mathbf S$ has the following behavior at the origin: as $z\rightarrow0$,  $z\in\mathbb{C}\backslash\Sigma$,
\[
\mathbf S(z)=
\begin{cases}
\mathcal O
\begin{pmatrix}
1 & \log\left\vert z\right\vert \\
1 & \log\left\vert z\right\vert
\end{pmatrix},
 & \text{as }z\rightarrow 0\text{ from the outer domain;}\\
\mathcal O
\begin{pmatrix}
\log\left\vert z\right\vert  & \log\left\vert z\right\vert \\
\log\left\vert z\right\vert  & \log\left\vert z\right\vert
\end{pmatrix},
  & \text{as }z\rightarrow 0\text{ from the inner domains.}%
\end{cases}
\]

\end{enumerate}

\subsection{The Szeg\H{o} function for $w_c$} \label{sec_Szego}

In this section we analyze in detail the structure and properties of the Szeg\H{o} function introduced in Section \ref{sec:intro}. Recall that for a non-negative function $h$ on $(-1,1)$ satisfying the Szeg\H{o} condition
$$
\int_{-1}^1  \frac{\log h(t)}{\sqrt{1-t^2}} \, dt>-\infty,
$$
we define in $\C\setminus [-1,1]$ its Szeg\H{o} function $D(\cdot, h)$ by
\begin{equation}\label{def_Szego}
    D(z,h)\isdef \exp \left(\frac{\sqrt{z^2-1}}{2\pi}\, \int_{-1}^1 \frac{\log h(t)}{\sqrt{1-t^2}}\, \frac{dt}{z-t}\right)= \exp\left( \sqrt{1-z^2}\,  \mathcal C   \left( \frac{\log h(t)}{  \sqrt{1-t^2}} \right) (z) \right),
\end{equation}
with $(\sqrt{1-z^2})_+>0$ for $z\in(-1,1)$ in the rightmost expression in \eqref{def_Szego}.

Due to the multiplicative property of the Szeg\H{o} function, we have that for $w_c$ defined in \eqref{defwanalytic},
\begin{equation}\label{SzegoTotal}
    D(z, w_c) = D(z, w) D(z, \Xi_c)\,.
\end{equation}
Straightforward computation shows that
\begin{equation}\label{SzegoPartial}
    D(z, w) = D(z, h)\, \dfrac{\left(  z-1\right)
^{ \alpha/2 }\left(  z+1\right)  ^{\beta/2}}{
\varphi^\frac{\alpha+\beta}{2}(z)   }, \quad D(z, \Xi_c) = c \, \exp\left(  -\lambda \log\left( \dfrac{1-i\sqrt{z^{2}-1}}{z}\right)\right)  ,
\end{equation}
where $D(\cdot, h)$ is computed by formula \eqref{def_Szego}, and
\begin{equation}\label{def_lambda}
 \lambda\isdef   i\dfrac{\log c}{\pi}\,.
\end{equation}
We must clarify that in \eqref{SzegoPartial} we take the main branches of $\left(  z-1\right)
^{ \alpha/2 }$, $\left(  z+1\right)  ^{\beta/2}$ and $\sqrt{z^{2}-1}$ that are positive for $z>1$, as well as the main branch of the logarithm.

From \eqref{SzegoTotal} we obtain that
\begin{equation}
D_{\infty}\isdef D(\infty, w_c)=\sqrt{c} \, D \left(  \infty, h \right)    2
^{-(\alpha+\beta)/2}>0. \label{Doo}%
\end{equation}

Let us study the boundary behavior of the Szeg\H{o} function on the interval. By \eqref{SzegoPartial},
$$
\lim_{\stackrel{z\to x \in (-1, 1),}{   \Im z>0}} D(z, w) =   e^{\pi i \alpha/2 } \left(  1-x\right)  ^{\alpha/2}\left(  1+x\right)  ^{\beta/2} \varphi_+^{-\frac{\alpha+\beta}{2} }(x) \lim_{\stackrel{z\to x \in (-1, 1),}{   \Im z>0}} D(z, h) ,
$$
where
\begin{equation}\label{phiplus}
    \varphi_+(x)=x+i\sqrt{1-x^2} = e^{  i\arccos (x) } ,
\end{equation}
with $\left(  1-x\right)  ^{\alpha/2}$, $\left(  1+x\right)  ^{\beta/2}$ and $\sqrt{1-x^2}$ positive for $x\in (-1,1)$.

Analogously,
$$
\lim_{\stackrel{z\to x \in (-1, 1),}{  \Im z<0}} D(z, w) =   e^{-\pi i \alpha/2 } \left(  1-x\right)  ^{\alpha/2}\left(  1+x\right)  ^{\beta/2} \varphi_+^{\frac{\alpha+\beta}{2} }(x) \lim_{\stackrel{z\to x \in (-1, 1),}{   \Im z<0}} D(z, h) .
$$

We can be more specific about the limit values of $D \left(  z, h\right) $ on $(-1,1)$ if we use the Sokhotskii-Plemelj formulas \cite[Section 4.2]{Gakhov:90}:
$$
\mathcal C _\pm \left( \frac{\log h(t)}{  \sqrt{1-t^2}} \right) (z) =\pm \frac{1}{2}\, \frac{\log
h\left(  t\right)  }{\sqrt{1-t^{2}}}  + \frac{1}{2 \pi i}\fint_{-1}^{1}\frac{\log h\left(  t\right)
}{\sqrt{1-t^{2}}}\frac{dt}{t-x},
$$
where $\fint$ is the integral understood in terms of its principal value. So, if we define $\hbar(x)$ as in \eqref{def_hbar1},
then using \eqref{def_Szego} we get
\begin{equation*}
\begin{split}
 \lim_{\stackrel{z\to x \in (-1, 1),}{  \Im z>0}}    D \left(  z, h\right)  & =  \sqrt{h(x)}\, e^{- i   \hbar(x)} ,  \\
  \lim_{\stackrel{z\to x \in (-1, 1),}{  \Im z<0}}   D \left(  z, h\right)  & =  \sqrt{h(x)}\, e^{i   \hbar(x)}.
\end{split}
\end{equation*}
Observe that $\hbar(x)$ is real-valued on $(-1,1)$, so that $\left|e^{\pm i  \hbar(x)} \right|=1$.
So, if we define on $(-1,1)$ the real-valued function
\begin{equation}\label{def_PhinoC}
    \Phi(x) \isdef \frac{\pi \alpha}{2} -\frac{\alpha +\beta }{2}\, \arccos(x) -  \hbar(x) ,
\end{equation}
then
\begin{equation*}
 \lim_{\stackrel{z\to x\in  (-1, 1), }{\pm \Im z>0}} D(z, w) = \sqrt{w(x)}\,   \exp\left(  \pm i \Phi(x) \right).
\end{equation*}

On the other hand, it is easy to check that  with the specified selection of the branch of the square root,
$$
z \mapsto \dfrac{1-i\sqrt{z^{2}-1}}{z}
$$
is a conformal mapping of $\C\setminus [-1,1]$ onto the lower half plane, such that the lower shore of $(-1,1)$ is mapped onto itself, while the upper boundary is mapped onto $(-\infty, -1)\cup (1, \infty)$. In particular,
$$
\lim_{\stackrel{z\to x\in  (0,1), }{ \Im z\neq 0}} \arg\left(  \dfrac{1-i\sqrt{z^{2}-1}}{z}\right)= 0,\quad \lim_{\stackrel{z\to x\in  (-1, 0), }{  \Im z\neq 0}} \arg\left(  \dfrac{1-i\sqrt{z^{2}-1}}{z}\right)= -\pi.
$$
Hence,
\begin{equation*}\label{prop fun Sze d+ 01}
 \lim_{\stackrel{z\to x\in  (0,1), }{\pm \Im z>0}} D(z, \Xi_c) = c \, \exp\left(  -\lambda \log\left| \dfrac{1\pm \sqrt{1-x^{2}}}{x}\right| \right)= c \, \exp\left(  \mp \lambda \log\left| \dfrac{1+ \sqrt{1-x^{2}}}{x}\right| \right),
\end{equation*}
with $\sqrt{1-x^2}>0$ on $(-1,1)$. Taking into account that $e^{-\lambda \pi i}=c$, we also get
\begin{equation*}\label{prop fun Sze d+ -10}
 \lim_{\stackrel{z\to x\in  (-1, 0), }{\pm \Im z>0}} D(z, \Xi_c) =   \exp\left(  -\lambda \log\left| \dfrac{1\pm \sqrt{1-x^{2}}}{x}\right| \right)=\exp\left(  \mp \lambda \log\left| \dfrac{1+ \sqrt{1-x^{2}}}{x}\right| \right).
\end{equation*}
Both identities can be summarized by
\begin{equation*}
 \lim_{\stackrel{z\to x\in  (-1, 0) \cup (0,1), }{\pm \Im z>0}} D(z, \Xi_c) = \sqrt{\Xi_c(x)}\,   \exp\left(  \mp  i\, \dfrac{\log c}{\pi} \,  \log\left| \dfrac{1+ \sqrt{1-x^{2}}}{x}\right| \right).
\end{equation*}

In order to clarify the local behavior of $D(z, \Xi_c)$ at the origin we observe that for $z\in\mathbb{C}\backslash\left(  -\infty
,1\right]  $ function $D(z, \Xi_c)$ coincides with $$
c \, \exp\left(  -\lambda \log(1-i\sqrt{z^{2}-1}) + \lambda \log(z)\right),
$$
if we take there the main branch of $\log(z)$, so that
$$
e^{-\lambda \log (z)}D\left(
z, \Xi_{c}\right)  =c \, \exp\left(  -\lambda   \log\left(
1-i\sqrt{z^{2}-1}\right)    \right).
$$
Since
$$
\lim_{\stackrel{z\to 0,}{ \Im z>0}} \log\left(
1-i\sqrt{z^{2}-1}\right)  =\log(2)\,,
$$
it yields
$$
D\left(  z, \Xi_c \right)  =c^{1+\frac{i}{\pi}\, \log\left(z/2 \right)}   \left(  1+o\left(
1\right)  \right),  \text{ as }z\rightarrow 0\text{, }   \Im %
z>0.
$$ The case $\Im z<0$ can be deduced using the symmetry of $D\left(
\cdot, w_{c}\right)$ with respect to $\R$.

We can summarize our findings in the following lemma:
\begin{lemma}
\label{Lem-D+} The Szeg\H{o} function $D(\cdot, w)$ for the weight $w$ defined in \eqref{defwanalytic} exhibits the following boundary behavior:
\begin{equation}\label{boundaryValuesFinal}
 \lim_{\stackrel{z\to x\in  (-1, 1), }{\pm \Im z>0}} D(z, w) = \sqrt{w(x)}\,   \exp\left(  \pm i \Phi(x) \right),
\end{equation}
with the notation introduced in \eqref{def_hbar1} and \eqref{def_PhinoC}.

Furthermore, for the step function $\Xi_c$,
\begin{equation*}
 \lim_{\stackrel{z\to x\in  (-1, 0) \cup (0,1), }{\pm \Im z>0}} D(z, \Xi_c) = \sqrt{\Xi_c(x)}\,   \exp\left(  \mp  i\, \dfrac{\log c}{\pi} \,  \log\left| \dfrac{1+ \sqrt{1-x^{2}}}{x}\right| \right),
\end{equation*}
and
\begin{equation}
D\left(  z, \Xi_c \right)  =c^{1\pm \frac{i}{\pi}\, \log\left(z/2 \right)}   \left(  1+o\left(
1\right)  \right),  \text{ as }z\rightarrow 0\text{, }  \pm \Im %
z>0. \label{D-0}%
\end{equation}
\end{lemma}
Obviously, the boundary behavior of the Szeg\H{o} function $D(\cdot, w_c)$ at $(-1,1)$ can be deduced from this Lemma and \eqref{SzegoTotal}.

\subsection{Outer parametrix}

Since $\left\vert \varphi\left(  z\right)  \right\vert >1$ for $z\in\mathbb{C}\backslash\left[  -1,1\right]  $, the matrix $\mathbf S$ introduced at the end of Subsection \ref{sec:First transformations} has jumps across each contour $\gamma_i$ that are exponentially close to $\mathbf  I$, as long as we stay away from the singularities $\pm 1$ and $0$. So, we can expect $\mathbf S$ to behave similarly to the $2\times 2$ solution $\mathbf N$ of the following RH problem in this region.
\begin{enumerate}
\item[(N1)] $\mathbf N$ is analytic in $\mathbb{C}\backslash\left[  -1,1\right]$;

\item[(N2)] $\mathbf N$ satisfies the following jump relations on  $\left(
-1,0\right)  \cup (0,1) $:%
\[
\mathbf N_{+}(x)=\mathbf N_{-}(x)
\begin{pmatrix}
0 & w_c\left(  x\right) \\
-w_c\left(  x\right)  ^{-1} & 0
\end{pmatrix};
\]
\item[(N3)] As $z\rightarrow\infty$,
\[
\mathbf N(z)=\mathbf {I}+\OO\left(  \frac{1}{z}\right) .
\]
\end{enumerate}
An explicit solution of this problem is well-known (see e.g.~\cite{MR2000g:47048} and \cite[Section 5]{MR2087231}) and can be built in terms of the Szeg\H{o} function $D(\cdot, w_c)$ and its value at infinity defined in \eqref{Doo}:
\begin{equation}
\mathbf N\left(  z\right)  \isdef D_{\infty}^{\sigma_{3}} \mathbf A(z)
D\left(  z, w_c\right)  ^{-\sigma_{3}}, \label{sol-N}%
\end{equation}
where
\begin{equation}
\label{defmatrixA}
\mathbf A(z)\isdef \begin{pmatrix}
A_{11} &  A_{12}\\
-A_{12} & A_{11} \end{pmatrix}= \begin{pmatrix}
\frac{a\left(  z\right)  +a^{-1}\left(  z\right)  }{2} & \frac{a\left(
z\right)  -a^{-1}\left(  z\right)  }{2i}\\
\frac{a\left(  z\right)  -a^{-1}\left(  z\right)  }{-2i} & \frac{a\left(
z\right)  +a^{-1}\left(  z\right)  }{2}%
\end{pmatrix}, \quad a\left(  z\right)  \isdef \tfrac{\left(  z-1\right)  ^{1/4}}{\left(
z+1\right)  ^{1/4}},
\end{equation}
and we take the principal branches in such a way that $a$ is analytic in
$\mathbb{C}\backslash\left[  -1,1\right]  $ with $a\left(  z\right)
\rightarrow1$ as $z\rightarrow\infty$. For future reference it is convenient to notice that an alternative expression for the entries of $\mathbf A$ can be obtained using that
\begin{equation}\label{AndPhi}
\begin{split}
  A_{11}(z) & =  \frac{a\left(  z\right)  +a^{-1}\left(  z\right)  }{2}=\frac{\varphi\left(
z\right)  ^{1/2}}{\sqrt{2}\left(  z^{2}-1\right)  ^{1/4}},
\\
A_{12}(z) & =\frac{a\left(  z\right)  -a^{-1}\left(  z\right)  }{2i}=\frac{i\varphi\left(
z\right)  ^{-1/2}}{\sqrt{2}\left(  z^{2}-1\right)  ^{1/4}}=\frac{i}{\varphi(z)}\, A_{11}(z) ,
\end{split}
\end{equation}
where we take again the main branches of the roots.

It is known (see \cite{MR2087231}) that $\mathbf N$ does not match the behavior of $\mathbf S$ at the endpoints of the interval $[-1,1]$, requiring a separate analysis there. Moreover, comparing the local condition (S5) for $\mathbf S$ with the behavior of $D(\cdot, w_c)$ at the origin (see \eqref{D-0}) we conclude that a local analysis will be needed also at $z=0$.

\subsection{Local parametrices at the endpoints of the interval}

We fix a $\delta\in (0, 1/8)$ and for each $\zeta \in \{-1, 1\}$ we consider the neighborhood $U_\zeta =\left\{  z\in \C:\, \left\vert z-\zeta  \right\vert
<\delta\right\}  $ such that $U_\zeta $ lies entirely in the domain $U$ of analyticity of $h$. We construct a $2\times 2$ matrix-valued
function $\mathbf P_{\zeta }$ in $U_\zeta \setminus \Sigma$ that exhibits the same jumps on $\Sigma\cap U_\zeta $ and the same local behavior at $z=\zeta $ as $\mathbf S$, and that matches the matrix $\mathbf N$ on the boundary $\partial U_\zeta $. Namely,
\begin{enumerate}
\item[(P$_\zeta $1)] $\mathbf P_{\zeta }$ is holomorphic in $U_{\zeta }\backslash\Sigma$ and continuous up to the boundary.

\item[(P$_\zeta $2)] $\mathbf P_{\zeta }$ satisfies the following jump relations:
\begin{align*}
\mathbf P_{\zeta +}(z)  &  =\mathbf P_{\zeta -}(z)\begin{pmatrix}
1 & 0\\
\frac{1}{w_c\left(  z\right)  }\, \varphi\left(  z\right) ^{-2n} & 1
\end{pmatrix}, \quad \text{for }z\in U_{\zeta}\cap\left( \bigcup_{i=1}^4\gamma_i\right) \setminus \{ \zeta \} ;\\
\mathbf P_{\zeta +}(x)  &  =\mathbf P_{\zeta -}(x)\begin{pmatrix}
0 & w_c \left(  x\right) \\
-\frac{1}{w_c\left(  x\right)  } & 0
\end{pmatrix}, \quad \text{for }x\in U_{\zeta }\cap\left( (-1,1) \right) .
\end{align*}

\item[(P$_\zeta $3)] As $n\rightarrow\infty$,
\[
\mathbf P_{\zeta}(z)\mathbf N^{-1}\left(  z\right)  =\mathbf {I}+\OO\left(  \frac{1}{n}\right)
\text{  uniformly\ for\ }z\in\partial U_{\zeta }\backslash\Sigma.
\]

\item[(P$_\zeta $4)] $\mathbf P_{\zeta }$ has the following behavior as $z\rightarrow \zeta $,  $z\in
U_{\zeta }\backslash\Sigma$:  with $s=\alpha $ if $\zeta =1$ and $s=\beta $ if $\zeta =-1$,
\begin{itemize}
 \item for $s<0$,
$$
\mathbf P_{\zeta }(z)=\OO \begin{pmatrix}
1 &  \left\vert z-\zeta \right\vert^s \\
1 &  \left\vert z-\zeta \right\vert^s
\end{pmatrix};
$$
\item for $s=0$,
$$
\mathbf P_{\zeta }(z)=\OO \begin{pmatrix}
\log \left\vert z-\zeta \right\vert &  \log \left\vert z-\zeta \right\vert  \\
\log \left\vert z-\zeta \right\vert &  \log \left\vert z-\zeta \right\vert
\end{pmatrix};
$$
\item for $s>0$,
\[
\mathbf P_{\zeta }(z)=
\begin{cases}
\OO \begin{pmatrix}
1 & 1 \\
1 & 1
\end{pmatrix},
  & \text{  as }z\rightarrow \zeta \text{ from the outer domain;}\\
\OO \begin{pmatrix}
 \left\vert z-\zeta \right\vert^{-s}  & 1 \\
 \left\vert z-\zeta \right\vert^{-s}  & 1
\end{pmatrix}, & \text{ as }z\rightarrow \zeta \text{ from the inner domain.}%
\end{cases}
\]
\end{itemize}

\end{enumerate}

We skip the details of construction of $\mathbf P_{\pm 1 }$, that can be found in \cite{MR2087231}.

\subsection{Local parametrix at the origin} \label{subsec:local0}

We fix a $\delta\in (0, 1/8)$ and consider the neighborhood $U_0=\left\{  z\in \C:\, \left\vert z\right\vert
<\delta\right\}  $ such that $U_0$ lies entirely in the domain $U$ of analyticity of $h$. We construct a $2\times 2$ matrix-valued function $\mathbf P_{0}$ in $U_0\setminus \Sigma$ that exhibits the same jumps on $\Sigma\cap U_0$ and the same local behavior at $z=0$ as $\mathbf S$, and that matches the matrix $\mathbf N$ on the boundary $\partial U_0$. Namely,
\begin{enumerate}
\item[(P$_0$1)] $\mathbf P_{0}$ is holomorphic in $U_{0}\backslash\Sigma$ and continuous up to the boundary.

\item[(P$_0$2)] $\mathbf P_{0}$ satisfies the following jump relations:
\begin{align*}
\mathbf P_{0+}(z)  &  =\mathbf P_{0-}(z)\begin{pmatrix}
1 & 0\\
\frac{1}{w_c\left(  z\right)  }\, \varphi\left(  z\right) ^{-2n} & 1
\end{pmatrix}, \quad \text{for }z\in U_{0}\cap\left( \bigcup_{i=1}^4\gamma_i\right) \setminus \{ 0 \} ;\\
\mathbf P_{0+}(x)  &  =\mathbf P_{0-}(x)\begin{pmatrix}
0 & w_c \left(  x\right) \\
-\frac{1}{w_c\left(  x\right)  } & 0
\end{pmatrix}, \quad \text{for }x\in U_{0}\cap\left( (-1,0)\cup (0,1) \right) .
\end{align*}

\item[(P$_0$3)] As $n\rightarrow\infty$,
\[
\mathbf P_{0}(z)\mathbf N^{-1}\left(  z\right)  =\mathbf {I}+\OO\left(  \frac{1}{n}\right)
\text{  uniformly\ for\ }z\in\partial U_{0}\backslash\Sigma\text{.}%
\]

\item[(P$_0$4)] $\mathbf P_{0}$ has the following behavior as $z\rightarrow0$,  $z\in
U_{0}\backslash\Sigma$:
\[
\mathbf P_{0}(z)=
\begin{cases}
\OO \begin{pmatrix}
1 & \log\left\vert z\right\vert \\
1 & \log\left\vert z\right\vert
\end{pmatrix},
  & \text{  as }z\rightarrow0\text{ from the outer domain;}\\
\OO \begin{pmatrix}
\log\left\vert z\right\vert  & \log\left\vert z\right\vert \\
\log\left\vert z\right\vert  & \log\left\vert z\right\vert
\end{pmatrix}, & \text{ as }z\rightarrow0\text{ from the inner domain.}%
\end{cases}
\]
\end{enumerate}

We build the solution of this problem in two steps. First we obtain a matrix $\mathbf P^{\left(  1\right)  }$ that satisfies conditions (P$_0$1, P$_0$2, P$_0$4), and after that, using an additional freedom in the construction, we take care of the matching condition (P$_0$3).

Let us define at this point an auxiliary function $W$ holomorphic  in $U\setminus \left( (-\infty, -1] \cup [1, +\infty)   \right)$ given by (see \eqref{defwanalytic})
\begin{equation}
W\left(  z\right)  \isdef \sqrt{c \, w(z) }\, , \quad \text{such that } W(x)>0 \text{ for } x\in (-1,1).
\label{Def-Wnew}%
\end{equation}
Then
\begin{equation}
W\left(  x\right)  =\begin{cases}
\sqrt{w_c\left(  x\right) c } , &  -1<x<0;\\
\sqrt{w_c\left(  x\right) c^{-1} }, &  0\leq x<1.%
\end{cases}
\label{rel W w}%
\end{equation}

We construct the matrix function $\mathbf P_{0}$ in the following form:
\begin{equation}
\mathbf P_{0}\left(  z\right)  =\mathbf E_{n}\left(  z\right) \mathbf   P^{\left(  1\right)  }\left(
z\right)  W\left(  z\right)  ^{-\sigma_{3}}\varphi\left(  z\right)
^{-n\sigma_{3}}, \label{Def-P}%
\end{equation}
where $\mathbf E_n$ is an analytic matrix-valued function in $U_0$ (to be determined). Matrix $\mathbf P^{\left(  1\right)  }$ is analytic in $U_0\setminus \Sigma$; using the properties of $W$ and $\varphi$ it is easy to show that
\begin{equation}\label{jumpsForP1_1}
\mathbf P_{+}^{\left(  1\right)  }\left(  x\right)  =\mathbf P_{-}^{\left(  1\right)  }\left(
x\right)
\begin{cases}
\begin{pmatrix}
0 & 1/c \\
-c & 0
\end{pmatrix}, & x\in (-\delta, 0), \\[4mm]
\begin{pmatrix}
0 & c \\
-1/c & 0
\end{pmatrix}, & x\in (0, \delta),
\end{cases}
\end{equation}
and
\begin{equation}\label{jumpsForP1_2}
\mathbf P_{+}^{\left(  1\right)  }\left(  z\right)  =\mathbf P_{-}^{\left(  1\right)  }\left(
z\right) \begin{cases}
\begin{pmatrix}
1 & 0 \\
c & 1
\end{pmatrix}, & z\in (\gamma_1 \cup \gamma_2)\cap U_0\setminus \{0\}, \\[4mm]
\begin{pmatrix}
1 & 0 \\
1/c & 1
\end{pmatrix}, & z\in (\gamma_3 \cup \gamma_4)\cap U_0\setminus \{0\}.
\end{cases}
\end{equation}
Taking into account that $W\left(  z\right)  =\OO\left(  1\right)  $ and $\varphi\left(  z\right)
=\OO\left(  1\right)  $ as $z\rightarrow0$, we conclude also from (P$_0$4) that
 $\mathbf P^{\left(  1\right)  }$ has the following behavior at the origin: as $z\rightarrow0$,  $z\in\mathbb{C}\backslash\Sigma$,
\begin{equation}\label{jumpsForP1_3}
\mathbf P^{\left(  1\right)  }(z)=
\begin{cases}
\OO \begin{pmatrix}
1 & \log|z| \\
1 & \log |z|
\end{pmatrix}, &  \text{ from the outer domain,}\\[4mm]
\OO \begin{pmatrix}
\log|z| & \log|z| \\
\log|z| & \log |z|
\end{pmatrix}, &  \text{ from the inner domain.}
\end{cases}
\end{equation}

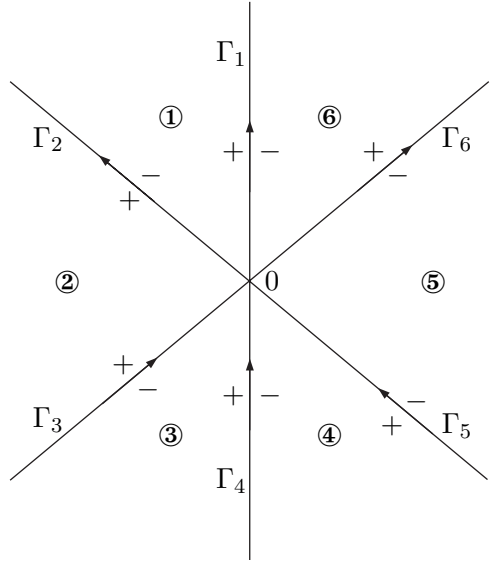
\begin{figure}[htb]
\centering \begin{overpic}[scale=1]%
{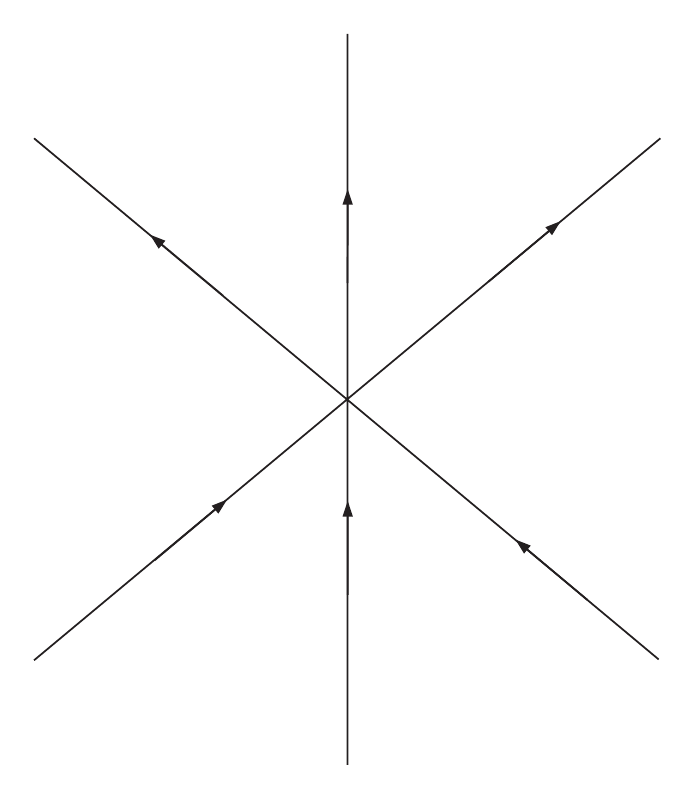}%
    \put(27,75){  \large \ding{172}  }
 \put(10,48){  \large  \ding{173}  }
 \put(27,23){  \large  \ding{174}  }
    \put(53,75){  \large \ding{177}  }
 \put(70,48){  \large  \ding{176}  }
 \put(53,23){  \large  \ding{175}  }
\put(21,35){$\small +$}
\put(25,31){$\small -$}
\put(62,70){$\small +$}
\put(66,66){$\small -$}
\put(39,70){$\small +$}
\put(45,70){$\small -$}
\put(39,30){$\small +$}
\put(45,30){$\small -$}
\put(22,62){$\small +$}
\put(25.5,66){$\small -$}
\put(65,25){$\small +$}
\put(69,29){$\small -$}
\put(46,48.5){$\small 0$}
\put(38,86){$\Gamma_1$}
\put(38,16){$\Gamma_4$}
\put(8,26){$\Gamma_3$}
\put(8,72){$\Gamma_2$}
\put(75,25){$\Gamma_5$}
\put(75,72){$\Gamma_6$}
\end{overpic}
\caption{Auxiliary contours $\Sigma_\psi$.}
\label{fig:psicontours}
\end{figure}

In order to construct $\mathbf P^{(1)}$ we solve first an auxiliary RH problem on a set $\Sigma_\Psi\isdef \bigcup_{j=1}^6 \Gamma_j$ of unbounded oriented straight lines converging at the origin, like in Fig.~\ref{fig:psicontours}.
More precisely,
\begin{align*}
 \Gamma_1 &=\left\{t e^{i\pi /2}:\, t>0 \right\}, \quad \Gamma_2 =\left\{t e^{3 i\pi /4}:\, t>0 \right\}, \quad \Gamma_3 =\left\{t e^{5 i\pi /4}:\, t>0 \right\},\\
 \Gamma_4 &=\left\{t e^{3 i\pi /2}:\, t>0 \right\}, \quad \Gamma_5 =\left\{t e^{- i\pi /4}:\, t>0 \right\}, \quad \Gamma_6 =\left\{t e^{  i\pi /4}:\, t>0 \right\}.
\end{align*}
These lines split the plane into 6 sectors, enumerated anti-clockwise from \ding{172} to \ding{177} as in Fig.~\ref{fig:psicontours}.

We look for a $2\times2$ matrix valued function
$\mathbf \Psi\left(  z\right)   $, satisfying the following conditions:
\begin{enumerate}
\item[($\Psi$1)] $\mathbf \Psi$ is analytic in $\mathbb{C}\backslash\Sigma_{\Psi}$.

\item[($\Psi$2)] for $k=1, \dots, 6$, $\mathbf \Psi$ satisfies the jump relation $\mathbf \Psi_{+}(\zeta)     =\mathbf \Psi_{-}(\zeta) J_k$ on $\Gamma_k$, with
\begin{equation}\label{Psi_k}
J_1=\begin{pmatrix}
0 & c\\
-1/c & 0
\end{pmatrix}, \quad J_2=J_6=\begin{pmatrix}
1 & 0\\
1/c & 1
\end{pmatrix}, \quad J_3=J_5=\begin{pmatrix}
1 & 0\\
c & 1
\end{pmatrix}, \quad J_4=\begin{pmatrix}
0 & 1/c\\
-c & 0
\end{pmatrix}.
\end{equation}

\item[($\Psi$3)] $\mathbf \Psi$ has the following behavior as $\zeta\rightarrow0$:%
\[
\mathbf \Psi(\zeta)=
\begin{cases}
\OO\begin{pmatrix}
\log| \zeta|  & \log| \zeta| \\
\log| \zeta|  & \log| \zeta|
\end{pmatrix},
 & \text{for }\zeta\in \text{\ding{172}}\cup \text{\ding{174}}\cup \text{\ding{175}}\cup \text{\ding{177}};\\[3mm]
\OO \begin{pmatrix}
1  & \log| \zeta| \\
1  & \log| \zeta|
\end{pmatrix},  & \text{as }\zeta\in \text{\ding{173}}\cup \text{\ding{176}}.
\end{cases}
\]
\end{enumerate}

If we use the notation $\lambda=i \log(c)/\pi$ introduced in \eqref{def_lambda}, then we readily see the connection of the RH problem above with that studied recently in \cite{Its07b}. Following the approach of \cite{Its07b} (with slight modifications), we construct $\mathbf \Psi$ explicitly in terms of the confluent hypergeometric functions
$$
\phi\left(  a,1;\zeta\right)\isdef {_{1}F_{1}}\left(  a;1;\zeta\right) \quad \text{and} \quad \psi\left(  a,1;\zeta\right) \isdef
\zeta^{-a}\,  {_{2}F_{0}}\left(  a, a; - ;-1/\zeta\right),
$$
that form a basis of solutions of the confluent hypergeometric equation $\zeta w^{\prime\prime}+\left(
1-\zeta \right)  w^{\prime}-aw=0$, see \cite[formula (13.1.1)]{abramowitz/stegun:1972}.
Namely, let
\begin{equation*}
  G\left(  a;\zeta\right)  \isdef  \phi\left(  a,1;\zeta\right) e^{-\zeta/2}, \quad H\left(  a;\zeta\right) \isdef \psi\left(  a,1;\zeta\right) e^{-\zeta/2} .
\end{equation*}
They are solutions of the confluent equation (see e.g.~\cite[formula (13.1.35)]{abramowitz/stegun:1972})
\begin{equation}
\zeta  w^{\prime\prime}+  w^{\prime}+ \left(  \frac{1}%
{2}-\frac{\zeta}{4}-a\right)  w=0; \label{def gen dif eq 1}%
\end{equation}
in fact, $G(a; \cdot)$ is the only entire solution of \eqref{def gen dif eq 1} such that $G(a; 0)=1$. Function $H(a,\zeta)$ is multivalued, and we take its principal branch in $-\frac{\pi}{2}<\arg\left(  \zeta\right)
<\frac{3\pi}{2}$. For these values of $\zeta$ we   define
$$
\widehat {\mathbf \Psi}\left(  \zeta\right)  \isdef
\begin{pmatrix}
\Gamma\left(  1-\lambda\right)  G\left(  \lambda;\zeta\right) & -H\left(  \lambda;\zeta\right) \\
\Gamma\left(  1+\lambda\right)   G\left( 1+ \lambda;\zeta\right) & \frac{\Gamma\left(  1+\lambda\right)  }{\Gamma\left(
-\lambda\right)  }\, H\left(  1+\lambda;\zeta\right)
\end{pmatrix}.
$$
By ($\Psi$2), if we set
\begin{equation*}
\mathbf \Psi\left(  \zeta\right)  \isdef
\begin{cases}
\mathbf {\widehat \Psi} \left(  \zeta\right)  J_6 J_1, & \text{for
}\zeta\in \text{\ding{172};} \\
\mathbf {\widehat \Psi}\left(  \zeta\right)  J_6 J_1 J_2 , & \text{for
}\zeta\in  \text{\ding{173};}\\
\mathbf {\widehat \Psi} \left(  \zeta\right)  J_6 J_1 J_2  J_3^{-1}, &
\text{for }\zeta\in \text{\ding{174};}\\
\mathbf {\widehat \Psi}\left(  \zeta\right)  J_5, & \text{for
}\zeta\in \text{\ding{175};}\\
\mathbf {\widehat \Psi} \left(  \zeta\right)  & \text{for }\zeta\in \text{\ding{176};}\\
\mathbf {\widehat \Psi} \left(  \zeta\right)  J_6 , & \text{for
}\zeta\in \text{\ding{177};}%
\end{cases}
\end{equation*}
then $\mathbf \Psi$ has the jumps across $\Sigma_\Psi$ specified in ($\Psi$2). Explicitly,
\begin{align}
\mathbf \Psi \left(  \zeta\right)  &=
\begin{pmatrix}
c^{-1}H\left(  \lambda;\zeta\right)  &
-\frac{\Gamma\left(  1-\lambda\right)  }{\Gamma\left(  \lambda\right)  }\,
H\left(  1-\lambda  ;e^{-\pi i}\zeta\right)  \\
-c^{-1}\frac{\Gamma\left(  1+\lambda\right)  }{\Gamma\left(
-\lambda\right)  }\, H\left(  1+\lambda;\zeta\right)   &
H\left(  -\lambda ;e^{-\pi i}\zeta\right)
\end{pmatrix}, \quad \zeta \in \text{\ding{172}},\label{sol Psi1}%
\\
\mathbf \Psi\left(  \zeta\right)  &=
\begin{pmatrix}
\Gamma\left(  1-\lambda\right)  G\left(  \lambda ;\zeta\right) & -\frac{\Gamma\left(  1-\lambda\right)  }{\Gamma\left(
\lambda\right)  }\, H\left(  1-\lambda ;e^{-\pi i}\zeta\right)   \\
\Gamma\left(  1+\lambda\right)  G\left(  1+\lambda ;\zeta\right)
  & H\left(  -\lambda,1;e^{-\pi i}\zeta\right)
\end{pmatrix}, \quad \zeta \in \text{\ding{173}},\label{sol Psi2}%
\\
\mathbf \Psi \left(  \zeta\right)  &=
\begin{pmatrix}
c H\left(  \lambda ;e^{-2\pi i}\zeta\right)   &
-\frac{\Gamma\left(  1-\lambda\right)  }{\Gamma\left(  \lambda\right)  }\,
H\left(  1-\lambda ;e^{-\pi i}\zeta\right)   \\
-c \frac{\Gamma\left(  1+\lambda\right)  }{\Gamma\left(  -\lambda
\right)  }\, H\left(  1+\lambda ;e^{-2\pi i}\zeta\right)   &
H \left(  -\lambda ;e^{-\pi i}\zeta\right)
\end{pmatrix}, \quad \zeta \in \text{\ding{174}},\label{sol Psi3}%
\\
\mathbf \Psi \left(  \zeta\right) & =
\begin{pmatrix}
-c \frac{\Gamma\left(  1-\lambda\right)  }{\Gamma\left(  \lambda\right)
}\, H\left(  1-\lambda ;e^{\pi i}\zeta\right)    & -H\left(
\lambda ;\zeta\right)   \\
c \, H\left(  -\lambda ;e^{\pi i}\zeta\right)    &
\frac{\Gamma\left(  1+\lambda\right)  }{\Gamma\left(  -\lambda\right)  }%
\, H\left(  1+\lambda;\zeta\right)
\end{pmatrix}, \quad \zeta \in \text{\ding{175}}, \label{sol Psi4}%
\\
\mathbf \Psi \left(  \zeta\right) & =
\begin{pmatrix}
\Gamma\left(  1-\lambda\right)  G \left(  \lambda ;\zeta\right)
  & -H \left(  \lambda ;\zeta\right)  \\
\Gamma\left(  1+\lambda\right)  G\left(  1+\lambda ;\zeta\right)
  & \frac{\Gamma\left(  1+\lambda\right)  }{\Gamma\left(
-\lambda\right)  }\, H\left(  1+\lambda ;\zeta\right)
\end{pmatrix}, \quad \zeta \in \text{\ding{176}}, \label{sol Psi5}%
\\
\mathbf \Psi\left(  \zeta\right) & =
\begin{pmatrix}
-c^{-1} \frac{\Gamma\left(  1-\lambda\right)  }{\Gamma\left(
\lambda\right)  }\, H\left(  1-\lambda ;e^{-\pi i}\zeta\right)    &
-H \left(  \lambda ;\zeta\right) \\
c^{-1} H\left(  -\lambda ;e^{-\pi i}\zeta\right)   &
\frac{\Gamma\left(  1+\lambda\right)  }{\Gamma\left(  -\lambda\right)  }%
\, H\left(  1+\lambda ;\zeta\right)
\end{pmatrix},\quad \zeta \in \text{\ding{177}}. \label{sol Psi6}%
\end{align}

Direct verification shows that $\mathbf \Psi$ coincides, after an appropriate change of parameters and a multiplication from the left by the constant matrix
\[
\begin{pmatrix}
1/c & 0\\
0 & 1
\end{pmatrix},
\]
with the solution of the corresponding RH problem found in \cite{Its07b} (cf.~formulas (7.26) and (7.27) therein).  In consequence, the matrix-valued function $\mathbf \Psi$ defined in \eqref{sol Psi1}--\eqref{sol Psi6} solves the RH problem ($\Psi$1)--($\Psi$3); moreover, $\det \mathbf \Psi \equiv 1$.

In order to construct the analytic function $\mathbf E_n$ in \eqref{Def-P} we need to study also the asymptotic behavior of $\mathbf \Psi$ at infinity. Let us introduce the notation
$$
\tau_\lambda  \isdef   \frac{ \Gamma(-\lambda ) }{\Gamma(\lambda )}.
$$
Then for purely imaginary values of $\lambda \neq 0$,
$$
\tau_\lambda  =   \frac{\overline{\Gamma(\lambda )}}{\Gamma(\lambda )}, \qquad |\tau_\lambda|=1, \qquad \tau_{-\lambda}=\overline{\tau_\lambda}.
$$
This value is not defined for $\lambda =0$; by continuity, we set $\tau _0=-1$.
\begin{lemma} \label{lemma_asymptoticsPsi}
As $\zeta\to \infty$, $\zeta\in \C\setminus \Sigma_\Psi$, and with the notation $\lambda=i \log(c)/\pi$, we have
\begin{equation}
\begin{split}
\mathbf \Psi\left(  \zeta\right)  =\left[   \mathbf  I  +\frac{\lambda}{\zeta}
\begin{pmatrix}
-\lambda  & -\tau_\lambda \\
-1/\tau_{\lambda}  & \lambda%
\end{pmatrix} +\OO\left(  \frac{1}{\left\vert \zeta\right\vert ^{2}}\right)  \right]
\zeta^{-\lambda\sigma_{3}}e^{  -\zeta\sigma_{3}/2} \\
\times \begin{cases}
 c^{-\sigma_{3}} , &\text{if } \frac{\pi}{2}<\arg\zeta<\frac{3\pi}{2}; \\
\begin{pmatrix}
0 & -1\\
1 & 0
\end{pmatrix}, & \text{if } -\frac{\pi}{2}<\arg\zeta<\frac{\pi}{2},
\end{cases}
\end{split}
 \label{Psir-ooNew}%
\end{equation}
where we use the main branch of $\zeta^{-\lambda}=e^{-\lambda \log \zeta}$ with the cut along $i\mathbb{R}_-$.
\end{lemma}
This result is a direct consequence of formulas (4.60)--(4.63) from \cite{Its07b}, and can be obtained by straightforward computation using the asymptotic properties of the confluent hypergeometric functions (see e.g.~\cite[formulas (13.5.1--2)]{abramowitz/stegun:1972}). In fact, formulas in \cite{abramowitz/stegun:1972} give us the complete expansion of $\mathbf \Psi$. 

Now we are ready to build $\mathbf P^{(1)}$ as in \eqref{Def-P}. Recall that $\varphi$ is a conformal mapping from $\C\setminus [-1,1]$ onto the exterior of the unit disk, so we can define in $\C\setminus \R$ the analytic function
\begin{equation}
f\left(  z\right)  \isdef \begin{cases}
\pi i -2 \log\varphi\left(  z\right)   ,  & \text{
for } \Im z>0, \\
\pi i + 2 \log\varphi\left(  z\right)  ,  & \text{
for } \Im z<0,
\end{cases}
 \label{Def-f}%
\end{equation}
where we take the main branch of the logarithm. Using that $\varphi_{+}\left(  x\right)  \varphi_{-}\left(  x\right)     =1$ on $( -1,1)$ we  conclude that $f_{+}\left(  x\right)     =f_{-}\left(  x\right)$ there, so that $f$ is holomorphic in $\C\setminus \left( (-\infty,-1] \cup [1, +\infty) \right)$. For $|z|<1$ we have
\begin{equation}
f\left(  z\right)  =2iz+\frac{1}{3}\, iz^{3}+\OO\left(  z^{5}\right)  ,\ \text{as
}z\rightarrow0\text{.} \label{f--0}%
\end{equation}
Hence, for $\delta>0$ sufficiently small, $f$ is a
conformal mapping of $U_{0}$. Moreover, by \eqref{phiplus},
\begin{equation}\label{fOnR}
    f(x)=2 i \arcsin(x), \quad x\in (-1,1),
\end{equation}
so that $f$ maps the real interval $\left(  -1,1\right)  $ one-to-one onto the purely imaginary interval
$\left(  -\pi i,\pi i\right)  $.

We can always deform our contours $\gamma_i$ close to $z=0$ in such a way that
$$
f\left( \gamma_1\cap U_0\right) \subset \Gamma_3, \quad f\left( \gamma_2\cap U_0\right) \subset \Gamma_5, \quad f\left( \gamma_3\cap U_0\right) \subset \Gamma_2, \quad f\left( \gamma_4\cap U_0\right) \subset \Gamma_6.
$$
With this convention, set
\begin{equation}
\mathbf P^{\left(  1\right)  }\left(  z\right)  \isdef \mathbf \Psi\left(  nf\left(  z\right)
\right), \quad z \in U_0.  \label{sol-P1}%
\end{equation}
By ($\Psi$1)--($\Psi$3) and \eqref{f--0}, this matrix-valued function has the jumps and the local behavior at $z=0$ specified in \eqref{jumpsForP1_1}--\eqref{jumpsForP1_3}. Taking into account the definition \eqref{Def-f} we get that
$$
e^{n f(z)} =e^{n \pi i} \, \varphi^{\mp 2n}(z), \quad \text{for } \pm \Im z >0.
$$
Hence, by Lemma \ref{lemma_asymptoticsPsi},
\begin{equation}
\begin{split}
\mathbf \Psi\left( nf\left(  z\right)\right)  =\left[   \mathbf  I  +\frac{\lambda }{n f(z)}
\begin{pmatrix}
-\lambda  & -\tau_\lambda \\
-1/\tau_{ \lambda}  & \lambda %
\end{pmatrix} +\OO\left(  \frac{1}{n ^{2}}\right)  \right]
(nf (z))^{-\lambda\sigma_{3}}i^{-n \sigma_3}   \\
\times \begin{cases}
c^{-\sigma_{3}}  \varphi^{n \sigma_3}(z)  , &\text{if } \Im z>0; \\
\begin{pmatrix}
0 & -1\\
1 & 0
\end{pmatrix} \varphi^{n \sigma_3}(z) , & \text{if } \Im z<0,
\end{cases}
\end{split}
 \label{AsymptoticsP1infty}%
\end{equation}
where the main branch of $\left[  nf\left(  z\right)  \right]  ^{\lambda}$ is taken with the cut along $\left(  -\infty,0\right]  $.
Since
\begin{equation*}
    \left[    f\left(  z\right)  \right]  ^{\lambda
}=\left\vert  f\left(  z\right)  \right\vert ^{\lambda}\exp\left(  -\frac{\log
c}{ \pi}\arg\left(   f\left(  z\right)  \right)  \right) ,
\end{equation*}
straightforward computations show that
\begin{equation}\label{jumpFonMinus}
\left[    f\left(  x\right)  \right] _\pm ^{\lambda
}=\begin{cases}
 \left\vert  f\left(  x\right)  \right\vert ^{\lambda} c^{-1/2} , & \text{for } 0<x<1,\\
\left\vert  f\left(  x\right)  \right\vert ^{\lambda}c^{-1/2 \mp 1} , & \text{for } -1<x<0,
\end{cases}
\end{equation}
where we assume the natural orientation of the interval.

Now we will build the analytic matrix $\mathbf E_{n}  $ in \eqref{Def-P}. In order to comply with condition (P$_0$3) above, we need
\[
\mathbf E_{n}\left(  z\right)  \sim \mathbf N\left(  z\right)  \varphi\left(  z\right)
^{n\sigma_{3}}W\left(  z\right)  ^{\sigma_{3}}\left[ \mathbf  P^{\left(  1\right)
}\left(  z\right)  \right]  ^{-1}%
\]
uniformly  for\ $z\in\partial U_{\delta }\backslash\Sigma$. Taking into account \eqref{AsymptoticsP1infty}, we define
\begin{equation}
\begin{split}
\mathbf E_{n}\left(  z\right) \isdef\mathbf  N\left(  z\right)  W\left(  z\right)  ^{\sigma_{3}}
\times \begin{cases}
i^{n\sigma_{3}}\left(nf\left(  z\right)  \right)  ^{\lambda\sigma_{3}} c^{ \sigma_{3}}  , &\text{if } \Im z>0; \\
i^{-n\sigma_{3}}\left(nf\left(  z\right)  \right)  ^{-\lambda\sigma_{3}} \begin{pmatrix}
0 & 1\\
-1 & 0
\end{pmatrix}, & \text{if } \Im z<0.
\end{cases}
\end{split}
 \label{DefinitionEN}%
\end{equation}
By construction, $\mathbf E_{n}$ is analytic in $U_0\backslash\mathbb{R}$. Furthermore, by (N2) and \eqref{rel W w}, for $x\in \left(
-\delta ,0\right)  \cup (0,\delta ) $,
\begin{equation*}
    \begin{split}
     W\left(  x\right)  ^{-\sigma_{3}}  \mathbf N_{-}^{-1}(x)\mathbf N_{+}(x) W\left(  x\right)  ^{\sigma_{3}} & =
\begin{pmatrix}
0 & w_c\left(  x\right)/W^2(x) \\
-W^2(x)/w_c\left(  x\right)   & 0
\end{pmatrix}\\ & =\begin{pmatrix}
0 & c^{\pm 1} \\
-c^{\mp 1}   & 0
\end{pmatrix}, \quad \text{for } \pm \Re x>0.
    \end{split}
\end{equation*}
From \eqref{jumpFonMinus} and \eqref{DefinitionEN} it follows that
$$
\mathbf E_{n-}^{-1}\left(  x\right) \mathbf E_{n+}\left(  x\right)=\mathbf  I,\quad  \text{for }  x\in \left(
-\delta ,0\right)  \cup (0,\delta ) .
$$
So, the origin is the only possible isolated singularity of $\mathbf E_n$ in $U_0$.
\begin{proposition} \label{propE0}
$$
\lim_{z\to 0} \mathbf E_n(z)= \frac{\sqrt{2}}{2}\, D_{\infty}^{\sigma_{3}} \begin{pmatrix}
1 & 1 \\ -1 & 1
\end{pmatrix} e^{i\eta_n \sigma_3},
$$
with $\eta_n$ introduced in \eqref{defOfEta}.
In particular, $\mathbf E_n$ is analytic in $U_0$.
\end{proposition}
\begin{proof}
Since $\mathbf E_n$ is analytic in a neighborhood of $0$ with an at most algebraic singularity there, it is sufficient to analyze its limit as $z\to 0$ from the upper half plane.
By \eqref{D-0} and \eqref{f--0},
\begin{align*}
\lim_{\stackrel{z\to 0}{\Im z>0}} D\left(  z, \Xi_c \right)
 f\left(  z\right)    ^{-\lambda} =\lim_{\stackrel{z\to 0}{\Im z>0}}  c^{1+ \frac{i}{\pi}\, \log\left(z/2 \right) - \frac{i}{\pi}\, \log\left( f(z) \right)}     =c^{3/2  } 4^{-\lambda }.
\end{align*}
On the other hand, by \eqref{boundaryValuesFinal} and \eqref{Def-Wnew} ,
\begin{align*}
\lim_{\stackrel{z\to 0}{\Im z>0}} D\left(  z, w \right)
 W\left(  z\right)    ^{-1} = c^{-1/2  } e^{i\Phi(0)}=  c^{-1/2  }\, \exp\left(i\frac{\alpha -\beta }{4}\, \pi - i\hbar(0) \right),
\end{align*}
with $\Phi$ given by \eqref{def_PhinoC} and $\hbar$ defined in \eqref{def_hbar1}.

Summarizing,
\begin{align*}
\lim_{\stackrel{z\to 0}{\Im z>0}} D\left(  z, w_c \right)^{-1}
 W\left(  z\right)  f\left(  z\right)    ^{ \lambda}  =\frac{4^\lambda}{c} \,  e^{-i\Phi(0)}.
\end{align*}
By \eqref{sol-N} and \eqref{DefinitionEN}, if $\Im z>0$,
\begin{equation}\label{expressionForE_n}
  \mathbf   E_n(z) = D_{\infty}^{\sigma_{3}}\,
\mathbf A(z) \, m_{n}(z)^{\sigma_3},
\end{equation}
with
\begin{equation}\label{def_m}
    m_n(z)\isdef \frac{W\left(  z\right)  f\left(  z\right)    ^{ \lambda}}{D\left(  z, w_c \right)} \, i^n n^\lambda c.
\end{equation}

Gathering the limits computed above, and using that
\begin{align*}
\lim_{\stackrel{z\to 0}{\Im z>0}} a(z)  = e^{\pi i/4}
\end{align*}
and the definition of $\eta_n$, the statement follows.
\end{proof}

Therefore, by construction the matrix-valued function $\mathbf P_0$ given by \eqref{Def-P} satisfies conditions (P$_{0}$1)--(P$_{0}$4). Moreover, it is easy to check that
\begin{equation*}
 \det \mathbf P_{0}\left( z\right)  =1 \quad \text{for every }z\in U_{0}\backslash\Sigma.
\end{equation*}

 \subsection{Final transformation}

\begin{figure}[htb]
\centering \begin{overpic}[scale=1.5]%
{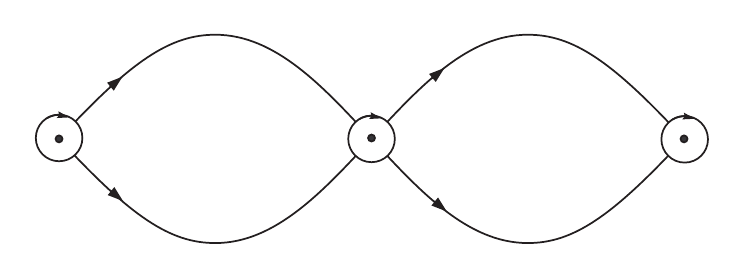}%
     \put(4,12){$\small -1 $}
\put(92,12){$\small 1$}
 \put(49,12){$\small 0 $}
\end{overpic}
\caption{Contours $\Sigma_R$.}
\label{fig:lenses3}
\end{figure}

Recall that matrices $\mathbf N$ and $\mathbf P_\zeta$, $\zeta \in \{-1, 0, 1 \}$ have $\det=1$ in their domains of definition. We may define
\begin{equation}
\mathbf R\left(  z\right)  \isdef
\begin{cases}
\mathbf S\left(  z\right)  \mathbf N^{-1}\left(  z\right),  &  z\in\mathbb{C}%
\backslash\left\{  \Sigma\cup U_{-1}\cup U_{0}\cup U_{1}\right\} ;\\
\mathbf S\left(  z\right) \mathbf  P_{\zeta}^{-1}\left(  z\right),  &  z\in
U_{\zeta}\setminus \Sigma, \; \zeta \in \{-1, 0, 1 \}.
\end{cases}
  \label{Def-R}%
\end{equation}
$\mathbf R$ is analytic in $\mathbb{C}\backslash\left\{  \Sigma\cup\partial U_{-1}\cup\partial
U_{0}\cup\partial U_{1}\right\}  $. In fact, since $\mathbf N$ matches the jump of $\mathbf S$ on $(-1,1)$, and $\mathbf P_\zeta$ matches the jumps of $\mathbf S$ within $U_\zeta$, $\zeta \in \{-1, 0, 1 \}$, we conclude that $\mathbf R$ is analytic in the complement to the contours $\Sigma_R$ depicted in Fig.~\ref{fig:lenses3}, with additional possible singularities at $\left\{  -1,0,1\right\}  $. But taking into account (S5) and the local behavior of $\mathbf P_ \zeta $ at these points (see (P$_ \zeta $4)), we conclude that these singularities are removable.

Now we compute the jumps of $\mathbf R$.  For the sake of brevity, we denote
\[
\Sigma_{R}^{out}\isdef \Sigma_{R}\setminus\left(  \partial U_{-1}%
\cup\partial U_{0}\cup\partial U_{1}\right) .
\]
Then by (S2) and \eqref{Def-R}, for $z\in\Sigma_{R}^{out}$,
\begin{equation}\label{jumpRelationsForR1}
 \mathbf  R_{+}\left(  z\right) =\mathbf R_{-}\left(  z\right) \mathbf  N\left(  z\right)
\begin{pmatrix}
1 & 0\\
 w_c\left(  z\right)^{-1}    \varphi\left(  z\right)^{-2n} & 1
\end{pmatrix}
\mathbf  N^{-1}\left(  z\right)  .
\end{equation}
On the other hand, for $\partial U_{j}$ ($j\in \{-1,0,1\}$) oriented clockwise, we have that $\mathbf R_{+}\left(  z\right)  =\mathbf S_{+}\left(  z\right) \mathbf  N^{-1}\left(  z\right)  $
and $\mathbf R_{-}\left(  z\right)  =\mathbf S_{-}\left(  z\right) \mathbf  P_{j}^{-1}\left(
z\right)  $. Hence,
\begin{equation}\label{jumpRelationsForR2}
\mathbf R_{+}\left(  z\right) =\mathbf R_{-}\left(  z\right) \mathbf  P_{j}\left(  z\right)  \mathbf N^{-1}\left(  z\right)  , \quad z\in
\partial U_{j},\quad j\in \{-1,0,1\}.
\end{equation}

Summarizing, $\mathbf R$ defined in \eqref{Def-R} is analytic in $\mathbb{C}\backslash\Sigma_{R}$,  satisfies the jump relations \eqref{jumpRelationsForR1}--\eqref{jumpRelationsForR2} on $\Sigma_{R}$, and has the following behavior as $z\to\infty$:%
\[
\mathbf R(z)=\mathbf {I}+\OO\left(  \frac{1}{z}\right) .
\]

By \eqref{jumpRelationsForR2} and (P$_0$3), as $n\to \infty$,
\begin{equation}
\mathbf R_{+}(z)=\mathbf R_{-}\left(  z\right)  \left(  \mathbf {I}+\OO\left(  \frac{1}%
{n}\right)  \right)  \quad \text{uniformly on } \partial U_{-1} \cup \partial U_{0} \cup \partial U_{1}.
\label{R n 00 2}%
\end{equation}
On the other hand, there exists a constant $0<q<1$ such that $\left\vert \varphi\left(  z\right)
\right\vert^{-1} \leq q <1$ uniformly on $\Sigma_{R}^{out}$. Since $\mathbf N$ does not depend on $n$, we conclude from \eqref{jumpRelationsForR1} that as $n\to \infty$,
\begin{equation}
\mathbf R_{+}(z)=\mathbf R_{-}\left(  z\right)  \left(  \mathbf {I}+\OO\left(  q^{2n}\right)
\right)  \quad \text{uniformly on }  \Sigma_{R}^{out}. \label{R n 00 1}%
\end{equation}
Motivated by \eqref{jumpRelationsForR1} --\eqref{R n 00 1} we define
\[
\mathbf \Delta\left(  s\right)  \isdef
\begin{cases}
\mathbf N\left(  s\right)
\begin{pmatrix}
1 & 0\\
 w_c\left(  s\right)^{-1}    \varphi\left(  s\right)^{-2n} & 1
\end{pmatrix}
\mathbf  N^{-1}\left(  s\right) - \mathbf{I},  & \text{for }s\in\Sigma_{R}^{out}; \\
\mathbf P_{\zeta}\left(  s\right) \mathbf  N^{-1}\left(  s\right) -\mathbf {I}, & \text{for }s\in\partial
U_{\zeta}, \; j\in \{-1,0,1 \},
\end{cases}
\]
so that $\mathbf R_{+}(z)=\mathbf R_{-}(z) (\mathbf  I + \mathbf  \Delta(z))$, $z\in \Sigma_{R}$.
Following \cite[Section 8]{MR2087231} we can show that $\mathbf \Delta$ has an asymptotic expansion in powers of $1/n$ of the form
\begin{equation}
\mathbf \Delta\left(  s\right)  \sim\sum_{k=1}^{\infty}\frac{\mathbf \Delta_{k}\left(
s,n\right)  }{n^{k}},\text{ \ \ as }n\rightarrow\infty,\text{ uniformly for
}s\in\Sigma_{R}. \label{DD-n-00}%
\end{equation}
By \eqref{R n 00 1}, for $k\in\mathbb{N}$,
\begin{equation}\label{delta=0}
\mathbf\Delta_{k}\left(  s\right)  =0\text{, \ \ for }s\in\Sigma_{R}^{out}\text{.}%
\end{equation}
Furthermore, by \cite[formulas (8.5)--(8.6)]{MR2087231},
\begin{align*}
\mathbf    \Delta_k(s) & =   \frac{(\alpha,k-1)}{2^k [\log\varphi\left(  s\right)]^k} \, \mathbf N(s) \left[e^{\pm \frac{ i \pi \alpha}{2}} c^{\frac{1}{2}} W(s)\right]^{\sigma_3}
    \begin{pmatrix}
        \frac{(-1)^k}{k} (\alpha^2 + \frac{1}{2}k-\frac{1}{4}) & -(k-\frac{1}{2})i \\[1ex]
        (-1)^k (k-\frac{1}{2}) i & \frac{1}{k} (\alpha^2 + \frac{1}{2}k - \frac{1}{4})
    \end{pmatrix} \\[1ex]
    &  \times   \left[e^{\pm \frac{i \pi \alpha}{2}} c^{\frac{1}{2}} W(s)\right]^{-\sigma_3} \mathbf N^{-1}(s),
    \qquad \mbox{for $\pm \Im s >0$ and $s\in\partial U_{1}$,}
\end{align*}
and
\begin{align*}
 \mathbf   \Delta_k(s) & =   \frac{(\beta,k-1)}{2^k [\log\left(
-\varphi\left(  s\right)\right) ]^k} \, \mathbf N(s)\,\left[e^{\mp  \frac{i \pi \beta}{2}} c^{- \frac{1}{2}} W(s)\right]^{\sigma_3}
    \begin{pmatrix}
        \frac{(-1)^k}{k} (\beta^2 + \frac{1}{2}k-\frac{1}{4}) & (k-\frac{1}{2})i \\[1ex]
        (-1)^{k+1} (k-\frac{1}{2}) i & \frac{1}{k} (\alpha^2 + \frac{1}{2}k - \frac{1}{4})
    \end{pmatrix} \\[1ex]
    & \times \left[e^{\mp  \frac{i \pi \beta}{2}} c^{- \frac{1}{2}} W(s)\right]^{-\sigma_3} \, \mathbf N^{-1}(s),
    \qquad \mbox{for $\pm \Im s >0$ and  $s\in\partial U_{-1}$,}
\end{align*}
where $(\alpha,0)\isdef 1$,
\[
    (\alpha,k) \isdef  \frac{(4\alpha^2-1)(4\alpha^2-9) \cdots (4\alpha^2 - (2k-1)^2)}{2^{2k} k!}.
\]
Each $\mathbf \Delta_k$ on the small contours encircling $\pm 1$ is independent of $n$ and possesses a meromorphic continuation to $U_{-1}$ and $U_{1}$ with the only pole at $\pm1$ of order at most $\left[  (k+1)/2\right]  $. However, unlike in the case analyzed in \cite{MR2087231}, the existence of a jump in the weight  is revealed through the contribution of the local parametrix $\mathbf P_{0}$, and hence, each $\mathbf \Delta_k$ is in general not independent on $n$, although uniformly bounded in $n$.

So, it remains to determine $\mathbf \Delta_k$ on $\partial U_0$. Here we calculate explicitly only the first term, $\mathbf \Delta_1$.

Using \eqref{sol-N}, \eqref{Def-Wnew}, \eqref{Def-P}, \eqref{Def-f}, \eqref{AsymptoticsP1infty} and \eqref{DefinitionEN},  we obtain
\begin{align*}
\mathbf  \Delta\left(  s\right)
&   =\mathbf  E_{n}(s)\left[  \frac{\lambda }{n f(s)}
\begin{pmatrix}
-\lambda  & -\tau_\lambda \\
-1/\tau_{ \lambda}  & \lambda %
\end{pmatrix} +\OO\left(  \frac{1}{n^{2}}\right)  \right]  \mathbf E_{n}^{-1}(s), \quad s\in\partial U_{0}, \quad n\rightarrow\infty.
\label{comp DD}%
\end{align*}
Let us define
\begin{equation}\label{DD-1}
  \mathbf \Delta_1\left(  s\right) \isdef \frac{\lambda}{  f(s)}\,  \mathbf E_{n}(s)
\begin{pmatrix}
-\lambda  & -\tau_\lambda \\
-1/\tau_{ \lambda}  & \lambda %
\end{pmatrix}  \mathbf   E_{n}^{-1}(s), \quad s\in\partial U_{0}.
\end{equation}
Using that by \eqref{DefinitionEN},
$$
\mathbf E_n(s)= \mathbf F(s) \,  \left(i^n n^\lambda  \right)  ^{ \sigma_{3}}= \mathbf F(s) \,  \left(i^n c^{\frac{i}{\pi}\, \log n}  \right)  ^{ \sigma_{3}},
$$
where
\begin{equation*} 
\begin{split}
\mathbf F\left(  s\right) \isdef  \begin{cases}
 \mathbf N\left(  s\right)  W\left(  s\right)  ^{\sigma_{3}} c^{ \sigma_{3}}   f\left( s\right)  ^{\lambda\sigma_{3}}  , &\text{if } \Im s>0; \\
\mathbf  N\left(  s\right)  W\left( s\right)  ^{\sigma_{3}}   \begin{pmatrix}
0 & 1\\
-1 & 0
\end{pmatrix}   f\left(  s\right)   ^{ \lambda\sigma_{3}}, & \text{if } \Im s<0,
\end{cases}
\end{split}
\end{equation*}
we conclude that for $s\in \partial U_0$,
\begin{equation}\label{delta1}
\mathbf \Delta_1\left(  s, n\right) = \frac{\lambda }{  f(s)}\,  \mathbf F (s)
\begin{pmatrix}
-\lambda  &  (-1)^{n+1} n^{2\lambda }\tau_\lambda \\
(-1)^{n +1} n^{-2\lambda } /\tau_{  \lambda}  & \lambda
\end{pmatrix}  \mathbf   F^{-1}(s),
\end{equation}
which is uniformly bounded in $n$, so that $\mathbf \Delta_{1}$ in \eqref{DD-1}--\eqref{delta1} is genuinely the first coefficient in the expansion \eqref{DD-n-00}.

Similar analysis can be performed for $\mathbf \Delta_k\left(  \cdot, n\right)$, $k\geq 2$, taking higher order terms in the expansion of $\mathbf \Psi$ in \eqref{Psir-ooNew}.

The explicit expression \eqref{delta1} and the local behavior of $f$ show that $\mathbf \Delta_{1}\left(  s, n\right)  $ has an analytic continuation to $U_{0}$ except for the origin, where it has a simple pole. Again, a similar conclusion is valid for other $\mathbf \Delta_{k}\left(  s, n\right)  $, except that now the pole is of order $k$.

As in \cite[Theorem 7.10]{MR2001f:42037} we obtain from \eqref{DD-n-00} that
\begin{equation}
\mathbf R\left(  z\right)  \sim\mathbf {I}+\sum_{j=1}^{\infty}\frac{\mathbf R^{\left(
j\right)  }\left(  z, n\right)  }{n^{j}},\text{ \ \ as }n\rightarrow
\infty\text{,} \label{R-n--00}%
\end{equation}
uniformly for $z\in\mathbb{C}\backslash\left\{  \partial U_{ -1}%
\cup\partial U_{ 0}\cup\partial U_{ 1}\right\}  $ where each
$\mathbf R^{\left(  j\right)  }\left(  z\right)  $ is analytic, uniformly bounded in $n$, and%
\[
\mathbf R^{\left(  j\right)  }\left(  z, n\right)  =\OO\left(  \frac{1}{z}\right)  \quad \text{as }z\rightarrow\infty\text{.}%
\]
This is a bona fide asymptotic expansion near infinity, since
\begin{equation*}
\forall l\in\mathbb{N}\ \exists C>0:\left\vert z\right\vert \geq
2\Rightarrow\left\Vert \mathbf R\left(  z\right)  -\mathbf{I}-\sum_{j=1}^{l}%
\frac{\mathbf R^{\left(  j\right)  }\left(  z, n\right)  }{n^{j}}\right\Vert \leq\frac
{C}{\left\vert z\right\vert n^{l+1}}\text{,} 
\end{equation*}
for any matrix norm $\left\Vert \cdot\right\Vert $. The proof is based on the integral representation for $\mathbf R$,
$$
\mathbf R(z)=\mathbf  I + \frac{1}{2\pi i}\, \int_{\Sigma_R} \frac{\mathbf R_-(s) \mathbf \Delta(s,n)}{s-z}\, ds, \quad z\in \C\setminus \Sigma_R
$$
(see \cite{MR2001f:42037}); although in our case the coefficients $\mathbf \Delta_k$ and $\mathbf R_k$ in \eqref{DD-n-00} and \eqref{R-n--00} depend on $n$, their uniform boundedness allows to follow the steps of the proof of Lemma 8.3 in \cite{MR2087231}.
In particular,  expanding the jump relation $\mathbf R_{+}=\mathbf R_{-}\left(  \mathbf {I}+\mathbf \Delta\right)
$ up to order $1/n$ we find that%
\begin{equation*}
\mathbf R_{+}^{\left(  1\right)  }\left(  s,n\right)  -\mathbf R_{-}^{\left(  1\right)
}\left(  s,n\right)  =\mathbf \Delta_{1}\left(  s,n\right)  , \quad \text{ for }s\in\partial
U_{ -1}\cup\partial U_{ 0}\cup\partial U_{ 1}.
\end{equation*}
Since $\mathbf R^{\left(  1\right)  }$ is analytic in the complement of $\partial
U_{-1}\cup\partial U_{0}\cup\partial U_{1}$ (see \eqref{delta=0}) and vanishes at infinity, by the Sokhotskii-Plemelj formulas,
\[
\mathbf R^{\left(  1\right)  }\left(  z,n\right)  =\frac{1}{2\pi i}\int_{\partial
U_{-1}\cup\partial U_{0}\cup\partial U_{1}}\frac
{\mathbf \Delta_{1}\left(  s, n\right)  }{s-z}\, ds.
\]
Recall that $\mathbf \Delta_{1}$ can be extended analytically inside $U_j$'s with simple poles at $\pm 1 $ and $0$; let us denote by $A^{\left(  1\right)  }\left(  n\right)  $,
$B^{\left(  1\right)  }\left(  n\right)  $ and $C^{\left(  1\right)  }\left(
n\right)  $ the residue of $\mathbf \Delta_{1}(\cdot, n)$ at $1$, $-1$ and $0$, respectively. Then residue calculus gives
\begin{equation}
\mathbf R^{\left(  1\right)  }\left(  z,n\right)  =
\begin{cases}
\dfrac{A^{\left(  1\right)  }\left(  n\right)  }{z-1}+\dfrac{B^{\left(
1\right)  }\left(  n\right)  }{z+1}+\dfrac{C^{\left(  1\right)  }\left(
n\right)  }{z}, & \text{for }z\in\mathbb{C}\backslash\left\{  U_{ -1}\cup
U_{ 0}\cup U_{ 1}\right\}  ;\\[3mm]
\dfrac{A^{\left(  1\right)  }\left(  n\right)  }{z-1}+\dfrac{B^{\left(
1\right)  }\left(  n\right)  }{z+1}+\dfrac{C^{\left(  1\right)  }\left(
n\right)  }{z}-\mathbf \Delta_{1}\left(  z,n\right) , & \text{for }z\in
U_{ -1}\cup U_{ 0}\cup U_{ 1}  .
\end{cases}
  \label{comp-R1}%
\end{equation}
Residues $A^{\left(  1\right)  }\left(  n\right)  $ and $B^{\left(  1\right)  }\left(  n\right)  $ are in fact independent of $n$; they have been determined in \cite[Section 8]{MR2087231}:
\begin{equation}
\begin{split}
\label{A-B}
A^{\left(  1\right)  }\left(  n\right) & =A^{(1)} =\frac{4\alpha^2-1}{16}\, D_{\infty}^{\sigma_3}
    \begin{pmatrix}
        -1 & i \\
        i & 1
    \end{pmatrix}D_{\infty}^{-\sigma_3},
\\
B^{\left(  1\right)  }\left(  n\right) & =B^{(1)} = 
    \frac{4\beta^2-1}{16}\, D_{\infty}^{\sigma_3}
    \begin{pmatrix}
        1 & i\\
        i & -1
    \end{pmatrix}
    D_{\infty}^{-\sigma_3}
\end{split}
\end{equation}
(notice however an extra factor $\sqrt{c}$ in the constant $D_\infty$ with respect to \cite{MR2087231}). The value of the remaining residue $C^{\left(
1\right)  }\left(  n\right)  $ is given in the following
\begin{proposition}
\label{prop coef C}
Coefficient $C^{\left(  1\right)  }\left(  n\right)  $ in \eqref{comp-R1} is given by
\[
C^{\left(  1\right)  }\left(  n\right)   = \frac{ \log c }{  2\pi} \, D_{\infty}^{\sigma_3} \begin{pmatrix}
-\cos  \theta_n   &  \lambda -i\sin  \theta_n      \\
 \lambda +i\sin  \theta_n      & \cos  \theta_n
\end{pmatrix}D_{\infty}^{-\sigma_3},
\]
where $ \theta_n$ is defined in \eqref{defOfTheta}.
\end{proposition}
\begin{proof}
Taking into account \eqref{f--0} and \eqref{DD-1} we conclude that
\begin{equation*}
   C^{\left(  1\right)  }\left(  n\right) = \frac{\lambda }{  2i}\,  \mathbf E_n (0)
\begin{pmatrix}
-\lambda  & -\tau_\lambda \\
-1/\tau_{\lambda}  & \lambda%
\end{pmatrix}   \mathbf E_n^{-1}(0).
\end{equation*}
By Proposition \ref{propE0},
\begin{align*}
 C^{\left(  1\right)  }\left(  n\right) & = \frac{ \lambda }{  4i}\,   D_{\infty}^{\sigma_{3}} \begin{pmatrix}
1 & 1 \\ -1 & 1
\end{pmatrix} \left(4^\lambda   e^{-i\Phi(0)} \, i^n n^\lambda  \right)^{\sigma_3} \begin{pmatrix}
-\lambda  & -\tau_\lambda \\
-1/\tau_{  \lambda}  & \lambda %
\end{pmatrix} \\ & \times   \left(4^\lambda   e^{-i\Phi(0)} \, i^n n^\lambda  \right)^{-\sigma_3}   \begin{pmatrix}
1 & -1 \\ 1 & 1
\end{pmatrix} D_{\infty}^{-\sigma_{3}} .
\end{align*}
With the notation \eqref{defOfEta} and choosing $\varsigma  \in \R$ such that $e^{i\varsigma}=\tau_\lambda $, we get
$$
\begin{pmatrix}
1 & 1 \\ -1 & 1
\end{pmatrix} e^{i \eta_n \sigma_3} \begin{pmatrix}
-\lambda  & -e^{i\varsigma} \\
-e^{-i\varsigma}  & \lambda
\end{pmatrix}    e^{-i \eta_n\sigma_3}   \begin{pmatrix}
1 & -1 \\ 1 & 1
\end{pmatrix} =2 \begin{pmatrix}
-\cos( 2\eta_n+\varsigma) & \lambda -i\sin(  2\eta_n+\varsigma) \\
\lambda +i\sin(  2\eta_n+\varsigma) & \cos(  2\eta_n+\varsigma)
\end{pmatrix}.
$$
It remains to observe that $2\eta _{n}+\varsigma=\theta_{n}$, and this settles the proof.
\end{proof}

\section{Asymptotic analysis. Proof of Theorems} \label{sec:proofs}

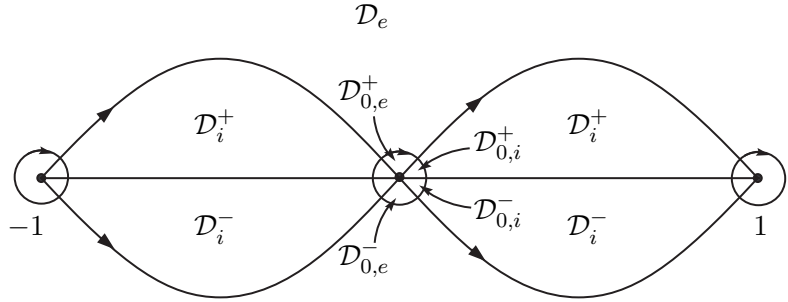
\begin{figure}[htb]
\centering \begin{overpic}[scale=1.5]%
{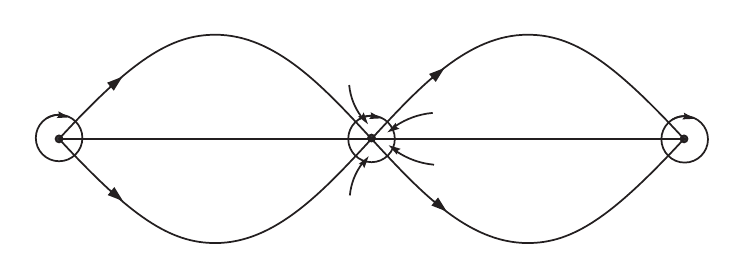}%
     \put(4,12){$\small -1 $}
\put(92,12){$\small 1$}
  \put(45,37){$\small \mathcal D_e$}
 \put(26,24){$\small \mathcal D_i^+$}
 \put(26,12){$\small \mathcal D_i^-$}
\put(70,24){$\small \mathcal D_i^+$}
 \put(70,12){$\small \mathcal D_i^-$}
\put(43,28){$\small \mathcal D_{0,e}^+$}
 \put(43,8.5){$\small \mathcal D_{0,e}^-$}
\put(59,22){$\small \mathcal D_{0,i}^+$}
 \put(59,14){$\small \mathcal D_{0,i}^-$}
\end{overpic}
\caption{Domains for $\mathbf Y$}
\label{fig:lenses4}
\end{figure}

Unraveling  the transformations $\mathbf Y\rightarrow \mathbf T\rightarrow \mathbf S\rightarrow \mathbf R$ we can obtain an expression for $\mathbf Y$.
We specify the following domains (see Fig.~\ref{fig:lenses4}):
\begin{itemize}
\item $\mathcal D_e$ is the unbounded component of $\mathbb{C}\backslash \Sigma_R$;
\item $\mathcal D_i^\pm$ correspond to the portion of the inner domain exterior to $U_ \zeta $, $\zeta \in \{-1,0,1\}$, lying in the upper (resp., lower) half-plane;
\item $\mathcal D_{\zeta ,e}^\pm$ is the subset of $U_ \zeta $ in the outer domain and upper (resp., lower) half plane;
 \item $\mathcal D_{\zeta,i}^\pm$ is the subset of $U_ \zeta $ in the inner domain and upper (resp., lower) half plane.
\end{itemize}

From (\ref{sol T}), (\ref{sol S}), and (\ref{Def-R}),
\begin{equation}
\mathbf Y\left(  z,n\right)  =
\begin{cases}
2^{-n\sigma_{3}}\mathbf R \mathbf N\varphi ^{n\sigma_{3}}\left(  z\right) , & z\in \mathcal D_e; \\
2^{-n\sigma_{3}}\mathbf R \mathbf N
\begin{pmatrix}
1 & 0\\
\pm\frac{1}{w_c}\, \varphi_{{}}^{-2n} & 1
\end{pmatrix}
 \varphi\left(  z\right)  ^{n\sigma_{3}}, & z\in \mathcal D_i; \\
2^{-n\sigma_{3}}\mathbf R \mathbf P_{\zeta}\varphi\left(  z\right)  ^{n\sigma_{3}} , & z\in \mathcal D_{\zeta,e}^\pm; \\
2^{-n\sigma_{3}}\mathbf R \mathbf P_{\zeta}
\begin{pmatrix}
1 & 0\\
\pm\frac{1}{w_c}\, \varphi_{{}}^{-2n} & 1
\end{pmatrix}
 \varphi\left(  z\right)  ^{n\sigma_{3}}, & z\in \mathcal D_{\zeta,i}^\pm;
\end{cases}
\label{sol Y final}%
\end{equation}
with $\zeta \in\left\{  -1,0,1\right\}  $.

Next, using the asymptotic expression for $\mathbf R$ derived above, we obtain information about the behavior of $\mathbf Y$ in different domains of the plane.

\subsection{Asymptotics for the monic orthogonal polynomials on $\mathbb{C}%
\backslash\lbrack-1,1]$. Proof of Theorem \ref{thm:outer}.} \label{subs:outerasympt}

If $K$ is a compact subset of $\mathcal D_e$, then by \eqref{R-n--00} and \eqref{sol Y final},
\begin{equation}\label{matrixYexteriorDomain}
  \mathbf   Y(z,n)=2^{-n\sigma_{3}}\mathbf R(z) \mathbf N(z)\varphi ^{n\sigma_{3}}\left(  z\right),
    \quad z \in K.
\end{equation}
Since $P_{n}\left(  z\right)  =\mathbf Y_{11}\left(
z,n\right)  $, we get by \eqref{sol-N}--\eqref{AndPhi} that
$$
\frac{2^{n}P_{n}\left(  z\right)  }{\varphi\left(  z\right)  ^{n}}=\frac{D_{\infty}}{D\left(  z, w_c \right)  } A_{11}(z) \mathfrak R(z),
$$
with
\begin{equation}\label{defRcal}
\mathfrak R(z)\isdef    \mathbf R  _{11}(z) -\frac{i}{D_{\infty}^{2} \varphi(z)}
\mathbf R  _{12}\left(  z\right).
\end{equation}
By \eqref{AndPhi}  and \eqref{R-n--00}, uniformly on $K$,
\[
\frac{2^{n}P_{n}\left(  z\right)  }{\varphi\left(  z\right)  ^{n}}=
\frac{D_{\infty}}{D\left(  z, w_c \right)  }\frac{\varphi\left(  z\right)  ^{1/2}%
}{\sqrt{2}\left(  z^{2}-1\right)  ^{1/4}}\left[  1+\frac{\mathcal R_n\left(  z\right)
}{n}+\OO\left(  \frac{1}{n^{2}}\right)  \right]  ,\text{ \ as }n\rightarrow
\infty\text{,}%
\]
with
\begin{equation}\label{def_H1}
 \mathcal R_n\left(  z\right)  \isdef \left( \mathbf  R^{\left(  1\right)  }\right)  _{11}\left(
z\right)  -\frac{i}{D_{\infty}^{2} \varphi(z)} \left(
\mathbf R^{\left(  1\right)  }\right)  _{12}\left(  z\right).
\end{equation}
Taking into account the expression for $\mathbf R^{(1)}$ in \eqref{comp-R1}, as well as \eqref{A-B} and Proposition \ref{prop coef C}, we get that
\begin{equation}\begin{split} \label{R112}
\mathbf  R^{(  1)  }  _{11}(z) & =   \frac{1-4\alpha^2 }{16(z-1)}  +\frac{4\beta^2-1}{16(z+1)} - \frac{ \log (c) \cos(\theta_n) }{  2\pi z}, \\
\mathbf  R^{(  1)  }  _{12}(z) & =  i D_\infty^2 \, \left(\frac{4\alpha^2-1 }{16(z-1) }  +  \frac{4\beta^2-1}{16(z+1) } + \frac{\log(c)}{2\pi }\, \frac{ \log (c)/\pi - \sin(\theta_n) }{   z}\right).
\end{split}
\end{equation}
The trivial identity $\varphi^2(z)+1=2z\,\varphi(z)$ yields
\begin{equation}
\label{identityForVarphi}
\frac{1}{z\pm 1}\, \left( 1\pm \frac{1}{\varphi(z)}\right)=\frac{2}{\varphi(z)\pm 1},
\end{equation}
and we conclude that in $K\subset \mathcal D_{e}$, $\mathcal R_n\left(  z\right) =\mathcal H_n\left(  z\right) $, with $\mathcal H_n  $ defined in \eqref{ValueofH}.

\subsection{Asymptotics of the recurrence coefficients} \label{sec:asymptRecurrencproofs}

Recall that monic polynomials $P_n$ satisfy the recurrence relation
\begin{equation*}
P_{n+1}(x)=\left(  x-b_{n}\right)  P_{n}(x)-a_{n}^{2}P_{n-1}(x),\quad n=0,1,\dots,
\end{equation*}
with $P_{-1}\left(  x\right)  =0$ and $a_n>0$. From \cite{Fokas92} (see also \cite{MR2000g:47048}) it follows that the coefficients can be found directly in terms of the matrix $\mathbf Y$ in \eqref{matrixYexteriorDomain}:
\begin{align}
a_{n}^{2}  &  =\lim\limits_{z\rightarrow\infty}z^{2}\mathbf Y_{12}\left(
z,n\right)  \mathbf Y_{21}\left(  z,n\right) \nonumber \\ & =\lim\limits_{z\rightarrow\infty}\left(  -\frac{D_{\infty}^{2}%
}{2i}+z\mathbf R_{12}\left(  z,n\right)  \right)  \left(  z\mathbf R_{21}\left(  z,n\right)
+\frac{1}{2iD_{\infty}^{2}}\right)  \text{,}\label{def an}\\
b_{n}  &  =\lim\limits_{z\rightarrow\infty}\left(  z-\mathbf Y_{11}\left(
z,n+1\right)  \mathbf Y_{22}\left(  z,n\right)  \right) =\lim\limits_{z\rightarrow\infty}z \left(  1-\mathbf  R_{11}\left(
z,n+1\right)  \mathbf R_{22}\left(  z,n\right)  \right)  \text{.} \label{def bn}%
\end{align}
We may take limits in the asymptotic expansion \eqref{R-n--00}; additionally to \eqref{R112} we have that
\begin{equation}
\begin{split}\label{R22}
\mathbf  R^{(  1)  }  _{21}(z) & =   \frac{i}{ D_\infty^{2}} \, \left(\frac{4\alpha^2-1 }{16(z-1) }  +  \frac{4\beta^2-1}{16(z+1) } + \frac{\log(c)}{2\pi }\, \frac{ \log (c)/\pi + \sin(\theta_n) }{   z}\right), \\
\mathbf  R^{(  1)  }  _{22}(z) & =  \frac{4\alpha^2 -1}{16(z-1)}  -\frac{4\beta^2-1}{16(z+1)} + \frac{ \log (c) \cos(\theta_n) }{  2\pi z}.
\end{split}
\end{equation}
Thus,
\begin{align*}
a_{n}^2  &  =\frac{1}{4}-\frac{\log c}{2 \pi n}\sin (\theta_{n})+\OO\left(  \frac{1}{n^{2}%
}\right),  \quad n\rightarrow\infty,
\end{align*}
which proves \eqref{asymptotics A}.
Analogously,
\begin{equation*}
\begin{split}
b_n & =\frac{\log c}{2 \pi }\, \frac{\cos (\theta_{n+1})- \cos (\theta_{n})}{n} +\OO\left(  \frac{1}{n^{2}%
}\right),
\quad n\rightarrow\infty.
\end{split}
\end{equation*}
By \eqref{defOfEta} and \eqref{defOfTheta},
$$
 \theta_{n+1}-\theta_n= \pi + 2 \frac{\log c}{\pi}\, \log\left( 1+\frac{1}{n}\right)   ,
$$
so that
\begin{align}
 b_n & =-\frac{\log c}{2 \pi }\, \frac{ \cos (\theta_n+  2 \frac{\log c}{\pi}\, \log( \frac{n+1}{n})  )+ \cos (\theta_{n})}{n} +\OO\left(  \frac{1}{n^{2}%
}\right) \nonumber \\
& = -\frac{\log c}{2 \pi n}\,   \left[ 2 \cos \left(\theta_n \right)  +\OO\left(  \frac{1}{n}\right) \right] +\OO\left(  \frac{1}{n^{2}}\right) ,\nonumber %
\end{align}
which proves \eqref{asymptotics B}.

In  \cite{Magnus1995} A.~Magnus conjectured that for the weight
\begin{equation*}
w\left(  x\right)  = \left(  1-x\right)  ^{\alpha} \left(  1+x\right)  ^{\beta}\left|
x_{0}-x\right|  ^{\gamma} \times
\begin{cases}
B, & \text{for }x\in\left[  -1,x_{0}\right)
\text{,}\\
A, & \text{for }x\in\left[  x_{0},1\right]  \text{,}%
\end{cases}
\end{equation*}
with $A$ and $B>0$ and $\alpha$, $\beta$ and $\gamma>-1$, and $x_{0}\in\left(
-1,1\right)  $, the recurrence coefficients of the corresponding orthogonal polynomials
exhibit the following behavior $n\rightarrow\infty$:
\begin{align}
a_{n}  &  =\frac{1}{2}-\frac{M}{n}\, \cos\left(  2n t_0-2\mu\log\left(
4n\sin t_0 \right)  -\widetilde \Phi\right)  +o\left(  1/n\right)  , \label{MagnusConj1}\\
b_{n}  &  =-\frac{2M}{n}\cos\left(   \left( 2 n+1\right)  t_0
-2\mu\log\left(  4n\sin t_0 \right)  -\widetilde  \Phi\right)  +o\left(
1/n\right)  ,\label{MagnusConj2}
\end{align}
where
\begin{align*}
x_0 & = \cos (t_0), \; 0<t_0<\pi, \quad \mu  =\frac{1}{2\pi}\log\frac{B}{A},\quad M=\frac{1}{2}%
\sqrt{ \frac{\gamma ^{2}}{4}  +\mu^{2}}\, \sin t_0, \\
\widetilde \Phi &  =\left(  \alpha+\frac{\gamma}{2}\right)  \pi-\left(  \alpha
+\beta+\gamma\right)  t_0 -2\arg\Gamma\left(  \frac{\gamma}{2}%
+i\mu\right)  -\arg\left(  \frac{\gamma}{2}+i\mu\right)  .
\end{align*}
Taking $B=1$, $A=c^2$, $\gamma =0$, and $x_0=0$ ($t_0=\pi/2$), we get $\mu=-\frac{ \log c}{ \pi}=i \lambda $, $M=|\log c|/(2\pi)$, and
\begin{align*}
\widetilde \Phi 
&  = \frac{\left(\alpha-  \beta \right)  \pi}{2}-2\arg\Gamma\left(
- \lambda \right)  -\arg\left(  - \lambda \right)\\
& =  \frac{\left(  \alpha-  \beta\right)  \pi}{2}+ 2\arg\Gamma\left(
 \lambda \right)  +\frac{\pi}{2}\, \sgn(\log (c)).
\end{align*}
Replacing these expressions in \eqref{MagnusConj1} and using the definition in \eqref{defOfEta} we obtain
\begin{align}
 a_n & =\frac{1}{2}-\frac{|\log c|}{2\pi n}\, \cos\left( \theta_n -\frac{\pi}{2}\, \sgn(\log (c))\right)  +o\left(  1/n\right) \nonumber \\
 & =\frac{1}{2}-\frac{ \log c }{2\pi n}\, \sin\left( \theta_n  \right)  +o\left(  1/n\right),   \label{asymptotics A Magnus}
\end{align}
and
\begin{align}
 b_n & = -\frac{|\log c|}{\pi n}\sin\left(   \theta_n    -\frac{\pi}{2}\, \sgn(\log (c)) \right)  +o\left(  1/n\right) \nonumber \\
 & =-\frac{\log c}{\pi n}\cos\left(   \theta_n    \right) +o\left(  1/n\right) \,.  \label{asymptotics B Magnus}
\end{align}
Comparing these expressions with \eqref{asymptotics A}--\eqref{asymptotics B} we see that Magnus' conjecture is valid for $\gamma=0$; moreover, we have shown that in this situation we can replace the error term $o(1/n)$ by a more precise $\OO(1/n^2)$.

\subsection{Asymptotics for the the leading coefficient $k_{n}$}

By \eqref{sol Y},
\[
k_{n}^{2}=-\frac{1}{2\pi i}\lim\limits_{z\rightarrow\infty}z^{-n}\mathbf Y_{21}\left(
z,n+1\right),
\]
and with \eqref{sol Y final},%
\[
k_{n}^{2}=-\frac{1}{2\pi i}\lim\limits_{z\rightarrow\infty}\left[
 \left(  \frac{2 \, \varphi\left(  z\right)  }{z}\right)  ^{n+1}\left(
z\mathbf R_{21}\left(  z,n+1\right) \mathbf  N_{11}\left(  z\right)  +z \mathbf R_{22}\left(
z,n+1\right) \mathbf  N_{21}\left(  z\right)  \right)  \right].
\]
Taking into account (N3), \eqref{AndPhi} and \eqref{R22}, we see that
\begin{align*}
&\lim_{z\to \infty} z \mathbf R_{21}\left(  z,n+1\right)  = \frac{i}{ n D_\infty^{2}} \, \left(\frac{2 \alpha^2+ 2\beta^2 -1 }{8  }    + \frac{\log(c)}{2\pi }\, \left( \frac{\log(c)}{ \pi } + \sin(\theta_{n+1})  \right) \right) + \OO\left( \frac{1}{n^2}\right),\\
&\lim_{z\to \infty} z \mathbf N_{21}\left(  z\right)  =-\frac{i}{2 D_\infty^{2}},
\end{align*}
and
\[
k_{n}^{2}= \frac{4^{n}}{ \pi D_\infty^{2} }   \left[
1  -  \left(\frac{2 \alpha^2+ 2\beta^2 -1 }{4  }    + \frac{\log(c)}{ \pi }\, \left( \frac{\log(c)}{ \pi } + \sin(\theta_{n+1})  \right)  \right)\,  \frac{1}{  n }   + \OO\left( \frac{1}{n^2}\right)  \right],
\]
and this proves Theorem \ref{thm:leading}.

\subsection{Asymptotics for the monic orthogonal polynomials in $U_0$ and on $\left(
-\delta,\delta\right).  $} \label{sec:asymptoninterval}

By analyticity of $P_n$'s, it is sufficient to consider $z\in \mathcal D_{0,i}^+$ and $\Re z>0$. Using formulas  \eqref{Def-P},  \eqref{sol-P1} and \eqref{sol Y final} we get
\begin{equation}
\mathbf Y\left(  z,n\right)  =2^{-n\sigma_{3}}\mathbf R\left(  z\right)\mathbf   E_{n}\left(
z\right) \mathbf  \Psi\left(  nf\left(  z\right)  \right)  W\left(  z\right)
^{-\sigma_{3}}\varphi\left(  z\right)  ^{-n\sigma_{3}}
\begin{pmatrix}
1 & 0\\
\frac{1}{w_c}\, \varphi^{-2n} & 1
\end{pmatrix}
  \varphi\left(  z\right)  ^{n\sigma_{3}}\text{.} \label{dem ass k 2}%
\end{equation}
We are interested in the first column of $\mathbf Y$, which is obtained multiplying the r.h.s.~of \eqref{dem ass k 2} from the right by the column vector $\left(  1, 0\right)^T  $. Observe that
\begin{equation*}
\begin{split}
W\left(  z\right)  ^{-\sigma_{3}}\varphi\left(  z\right)  ^{-n\sigma_{3}%
} \begin{pmatrix}
1 & 0\\
\frac{1}{w_c(z)}\, \varphi(z)^{-2n}  & 1
\end{pmatrix}
   \varphi\left(  z\right)  ^{n\sigma_{3}}
\begin{pmatrix}
1\\
0
\end{pmatrix} & =\begin{pmatrix}
1/W(z)\\
W(z)/w_c(z)
\end{pmatrix}\\ & =\frac{1}{W(z)}\, \begin{pmatrix}
1\\
1/c
\end{pmatrix}, 
\end{split}
\end{equation*}%
where we have taken into account the definition of $W$ in $\mathcal D_{0,i}^+$. Thus,
\begin{equation}\label{intermediate_local1}
    W(z) \mathbf Y\left(  z,n\right) \begin{pmatrix}
1\\
0
\end{pmatrix} =2^{-n\sigma_{3}}\mathbf R\left(  z\right)\mathbf   E_{n}\left(
z\right) \mathbf  \Psi\left(  nf\left(  z\right)  \right) \, \begin{pmatrix}
1\\
1/c
\end{pmatrix}.
\end{equation}
Notice that   $ \mathcal D_{0,i}^+$ is mapped by $f$ onto the sector denoted by \ding{172} in Figure \ref{fig:psicontours}, and vector $(1, 1/c)^T$ corresponds to the first column of the jump matrix $J_2$ in \eqref{Psi_k}. Taking into account ($\Psi2$) we conclude that the product of the last two matrices in the right hand side of \eqref{intermediate_local1} is equal to the first column of $\mathbf \Psi$ in \eqref{sol Psi2}:
\begin{equation}
\mathbf \Psi\left(  nf\left(  z\right)  \right) \, \begin{pmatrix}
1\\
1/c
\end{pmatrix} = \begin{pmatrix}
\Gamma\left(  1-\lambda\right)  G\left(  \lambda ;nf\left(  z\right)\right)   \\
\Gamma\left(  1+\lambda\right)  G\left(  1+\lambda ;nf\left(  z\right)\right)
\end{pmatrix}. \label{dem ass k 5}%
\end{equation}
By \eqref{expressionForE_n},
\begin{equation}\label{intermediateY}
  W(z) \mathbf Y\left(  z,n\right) \begin{pmatrix}
1\\
0
\end{pmatrix} =2^{-n\sigma_{3}}\mathbf R\left(  z\right)   D_{\infty}^{\sigma_{3}} \mathbf A(z)\,
 m_n(z)^{\sigma_3}\begin{pmatrix}
\Gamma\left(  1-\lambda\right)  G\left(  \lambda ;nf\left(  z\right)\right)   \\
\Gamma\left(  1+\lambda\right)  G\left(  1+\lambda ;nf\left(  z\right)\right)
\end{pmatrix},
\end{equation}
with $\mathbf A$ and $m_{n}$ defined in \eqref{defmatrixA} and \eqref{def_m}, respectively.
Taking into account formulas \eqref{AndPhi}, we conclude that
\begin{equation*}
\begin{split}
2^n W(z) P_n(z) & = D_\infty \, \frac{\varphi\left(
z\right)  ^{1/2}}{\sqrt{2}\left(  z^{2}-1\right)  ^{1/4}} \\
& \times \left\{ \mathfrak R(z) \, m_n(z) \Gamma\left(  1-\lambda\right)  G\left(  \lambda ;nf\left(  z\right)\right) + \widetilde{\mathfrak R}(z) \, m_n(z)^{-1} \Gamma\left(  1+\lambda\right)  G\left( 1+ \lambda ;nf\left(  z\right)\right) \right\},
\end{split}
\end{equation*}
where we have used notation \eqref{defRcal} and
$$
\widetilde{\mathfrak R}(z)\isdef \mathbf R_{11}(z)\frac{i}{  \varphi(z)} + \frac{1}{D_\infty^2  }\, \mathbf R_{12}(z).
$$

Inserting again  \eqref{R-n--00} we obtain the asymptotic expansion valid uniformly on compact subsets of $U_0$.  Using the function $\mathcal R_n$ defined in \eqref{def_H1} and introducing
\begin{equation*}
\begin{split}
    \widetilde {\mathcal R}_n\left(  z\right)  & \isdef \left( \mathbf  R^{\left(  1\right)  }\right)  _{11}\left(
z\right)  -\frac{i \varphi(z)}{D_{\infty}^{2} } \left(
\mathbf R^{\left(  1\right)  }\right)  _{12}\left(  z\right) 
\end{split}
\end{equation*}
we rewrite this identity for $P_n$  as
\begin{equation}\label{expressioncloseto0}
\begin{split}
 2^n P_n(z) W(z) =  & D_{\infty} A_{11}(z)
        \bigg[ \left( 1+\frac{1}{n}\, \mathcal R_n(z) +\OO\left(\frac{1}{n^{2}}\right)\right)  \, m_n(z) \Gamma\left(  1-\lambda\right)  G\left(  \lambda ;nf\left(  z\right)\right) \\
        & \left. + \frac{i}{\varphi(z)}\, \left( 1+\frac{1}{n}\,  \widetilde {\mathcal R}_n(z) +\OO\left(\frac{1}{n^{2}}\right)\right)   \,  \frac{\Gamma\left(  1+\lambda\right)  G\left( 1+ \lambda ;nf\left(  z\right)\right)}{m_n(z)}
        \right].
\end{split}
\end{equation}
Let us simplify this expression for the case when $z$ is on the real line. Taking the limit $z\to x\in (-\delta, \delta)$ from the upper half plane, we get by \eqref{defOfEta}, \eqref{fOnR}, \eqref{jumpFonMinus} and Lemma \ref{Lem-D+},
$$
   m_n(x)=      e^{ i \,  (\rho (x)+\eta_n)} , \quad  \text{for } x\in (-\delta,  \delta),
$$
with $\rho(x)$ given in \eqref{defFuncRho}, so that on $(-\delta, \delta)$,
$$
\overline{m_n(x)} = m_n(x)^{-1}.
$$
Additionally, we have $\overline{\lambda}=-\lambda $ and for $x\in \R$, on account of formulas (6.1.23) and (13.1.27) from \cite{abramowitz/stegun:1972}, respectively,
\begin{equation}\label{conjugates}
    \Gamma\left(  1+\lambda\right)   = \overline{\Gamma\left(  1-\lambda\right)    },\quad \text{and} \quad  G\left( 1+ \lambda ;ix \right)={_{1}F_{1}}\left(  -\lambda ;1;  -i x\right) e^{ix/2}= \overline{  G\left(  \lambda ;ix \right) }.
\end{equation}
Finally, on $(-\delta, \delta)$,
$$
\left( A_{11}(x)\right)_{+}= \frac{  \varphi_+(x)^{1/2}}{\sqrt{2} (x^2-1)_+^{1/4}}=\overline{\left( \frac{ i \varphi_+(x)^{-1/2}}{ \sqrt{2} (x^2-1)_+^{1/4}} \right)} = \overline{\left( A_{12}(x)\right)_{+}}=\frac{  e^{-i\arcsin(x)/2}}{  \sqrt{2} (1-x^2)^{1/4}}
$$
and by \eqref{delta1}, \eqref{comp-R1}, \eqref{A-B} and Proposition \ref{prop coef C},
$$
\overline{ \left( \mathbf  R^{\left(  1\right)  }\right)  _{11}\left(
x\right) } = \left( \mathbf  R^{\left(  1\right)  }\left(
x\right) \right)  _{11} , \quad  \overline{\left( \mathbf  R^{\left(  1\right)  }\right)  _{12}\left(
x\right) }=-\left( \mathbf  R^{\left(  1\right)  }\left(
x\right) \right)  _{12}.
$$
Gathering all this information in \eqref{expressioncloseto0} we conclude that locally uniformly on $(-\delta, \delta)$, as $n\to\infty$,
\begin{equation*}
\begin{split}
2^n P_n(x) W(x) &=   \frac{ \sqrt{2}\, D_{\infty} }{(1-x^2)^{1/4}}\\ & \times \Re \left[ \left( 1+ \frac{\mathcal R_{n}(x)}{n}  +\OO\left( \frac{1}{n^2}\right)\right) \, e^{-i\arcsin(x)/2} m_n(x) \Gamma\left(  1-\lambda\right)  G\left(  \lambda ;nf\left(  x\right)\right)\right],
\end{split}
\end{equation*}
with $\mathcal R_{n}$ defined in \eqref{def_H1}. Observe however that now the explicit expression for $\mathcal R_{n}$ differs from $\mathcal H_{n}$  defined in \eqref{ValueofH}: by \eqref{comp-R1}, in $\mathcal D_{0,i}^{+}$,
$$
\mathcal R_{n}(z)=\mathcal H_{n}- \left( \left( \mathbf  \Delta_{1}\left(
z,n\right)  \right)  _{11}-\frac{i}{D_{\infty}^{2}\varphi\left(  z\right)
}\left( \mathbf  \Delta_{1}\left(  z,n\right)  \right)  _{12}\right),
$$
with $\mathbf \Delta_{1}\left(  z,n\right)  $ given in \eqref{DD-1}. Using \eqref{AndPhi},
\eqref{expressionForE_n}, and \eqref{def_m}, we get
\begin{equation*}%
\begin{split}
\left(\mathbf   \Delta_{1}\left(  z,n\right)  \right)  _{11}  & =\frac{\lambda
}{f\left(  z\right)  }\frac{\varphi\left(  z\right)  ^{2}}{\varphi\left(
z\right)  ^{2}-1}\left(  -\lambda\left(  1+\frac{1}{\varphi\left(  z\right)
^{2}}\right)  -\frac{i}{\varphi\left(  z\right)  }\left(  \tau_{\lambda}%
m_{n}\left(  z\right)  ^{2}+\frac{1}{\tau_{\lambda}m_{n}\left(  z\right)
^{2}}\right)  \right) , \\
\left( \mathbf  \Delta_{1}\left(  z,n\right)  \right)  _{12}  & =\frac{\lambda
}{f\left(  z\right)  }\frac{\varphi\left(  z\right)  ^{2}}{\varphi\left(
z\right)  ^{2}-1}D_{\infty}^{2}\left(  2i\frac{\lambda}{\varphi\left(
z\right)  }-\tau_{\lambda}m_{n}\left(  z\right)  ^{2}-\frac{1}{\tau_{\lambda
}m_{n}\left(  z\right)  ^{2} \varphi\left(  z\right)  ^{2}}\right)  .
\end{split}
\end{equation*}%
Hence, by \eqref{identityForVarphi}  we have that
\begin{equation*}%
\begin{split}
 \mathcal{R} _{n}\left(  z\right)    & =\mathcal{H}_{n}\left(
z\right)  +\frac{\lambda}{f\left(  z\right)  }\left(  \lambda+\frac{i}%
{\varphi\left(  z\right)  }\overline{\tau_{\lambda}}m_{n}^{-2}\right)  .
\end{split}
\end{equation*}%
For further simplification of our formula we may take into account that by \cite[formula  (6.1.29)]{abramowitz/stegun:1972}, for $c\neq 1$,
\begin{align*}
\Gamma\left(  1-\lambda\right)   &=\overline{ \Gamma\left(  1+\lambda\right) } = \overline{ \lambda\Gamma\left(  \lambda\right)}
 =- i\frac{\log c}{\pi}\left\vert \Gamma\left(  \lambda\right)  \right\vert
e^{-i\arg\Gamma\left(  \lambda\right)  }\\
& =- i\frac{\log c}{\sqrt{\log c\sinh\left(  \log c\right)  }}e^{-i\arg
\Gamma\left(  \lambda\right)  }= - i \Upsilon(c)e^{-i\arg
\Gamma\left(  \lambda\right)  } , 
\end{align*}
with $\Upsilon(c)$ given in \eqref{defUpsilon},
and we obtain that for $x\in(-\delta,\delta)$,%
\begin{equation*}%
 \Gamma\left(  1-\lambda\right) m_{n}(x)   =  - i \Upsilon(c)e^{-i\arg
\Gamma\left(  \lambda\right)  }   e^{ i \,  (\rho (x)+\eta_n)}  =  - i \Upsilon(c)  e^{ i \,  (\rho (x)+\theta_{n}/2)}.
\end{equation*}%
Analogously,
\begin{equation*}%
\overline{\tau_{\lambda} } \, m_{n}(x)^{-2}=  e^{-2i\,(\rho(x)+\theta
_{n}/2)}, \quad x\in (-\delta, \delta),
\end{equation*}%
so that for $x\in (-\delta, \delta)$,
\begin{equation*}%
\begin{split}
\left(\frac{\lambda}{f\left(  z\right)  }\left(  \lambda+\frac{i}%
{\varphi\left(  z\right)  }\overline{\tau_{\lambda}}m_{n}^{-2}\right) \right)_{+ }&=
\frac{i \log c}{2\pi \arcsin(x)  }\left( \frac{\log c}{\pi }+ e^{-i\,(2\rho(x)+\theta
_{n}+ \arccos(x))}  \right).
\end{split}
\end{equation*}%
Summarizing,
\begin{equation*}
\begin{split}
2^n P_n(x) W(x) & = \frac{ \sqrt{2}\, D_{\infty}\Upsilon(c) }{(1-x^2)^{1/4}}\,\\
& \times  \Re \left[ \left( 1+ \frac{\mathcal R_{n}(x)}{n}  +\OO\left( \frac{1}{n^2}\right)\right) \,(-i) e^{-i\arcsin(x)/2}  e^{ i \,  (\rho (x)+\theta_{n}/2)}  G\left(  \lambda ;nf\left(  x\right)\right)\right],
\end{split}
\end{equation*}
which proves Theorem \ref{thm:local}.

Furthermore, with the appropriate rescaling and taking into account the local behavior of the terms in the right hand side of the asymptotic expression for $P_n$ we easily get the assertion of Corollary \ref{cor:localbeh}.

In order to prove Proposition \ref{prop:clockbehavior} we rewrite \eqref{asymptPicloseto0scaled} as
\begin{equation}\label{asymptPicloseto0scaledAlternative}
\begin{split}
P_n\left(\frac{\pi x}{n}\right) & =\frac{D_\infty\, \Upsilon(c)}{2^{n-1/2}\sqrt{c\, h(0)}}\, \left|  {_{1}F_{1}}\left(  \lambda;1;2\pi i x\right) \right|\, \Im \left[   e^{\frac{ i}{2}\, \left( \theta_n - \mathfrak G(2\pi x)\right) } \,     \left( 1 +   \OO\left( \frac{1}{n }\right)\right)  \right],
\end{split}
\end{equation}
where $\mathfrak G$ is the function introduced in \emph{(iii)} of Proposition \ref{prop:nozeros}, corresponding to $a=\log(c)/\pi$. Let us consider here only the case $c>1$ (the other case can be easily reduced to $c>1$ by a change of variables $x \mapsto -x$). Then $\mathfrak G$ is strictly increasing in $\R$.
If we denote by
\begin{equation}
\label{enumerationzeros}
\dots <\zeta_{-k}^{(n)} <\dots  <\zeta_{-1}^{(n)} <0\leq \zeta_{0}^{(n)}<\dots<\zeta_{k}^{(n)}<\dots
\end{equation}
the solutions of
$$
\frac{1}{2\pi}\,  \mathfrak G(2\pi x) \equiv \frac{\theta_n}{2\pi} \mod (\Z),
$$
then by \eqref{asymptPicloseto0scaledAlternative},
\begin{equation}
\label{xandy}
\lim_{n} \left(\frac{n }{\pi}\, x^{(n)}_{k} -  \zeta_{k}^{(n)} \right)=0, \quad k \in \Z,
\end{equation}
where we have used notation \eqref{enumerationzerosX}.
Since  $\mathfrak G(0)=0$, we have that $ \zeta_{0}^{(n)}$ is given by
$$
\frac{1}{2\pi}\, \mathfrak G(2\pi x) = \left\{ \frac{\theta_n}{2\pi}\right\},
$$
where $\{ \cdot\}$ is the fractional part of the number, which by strict monotonicity of $\mathfrak G$ shows that
\begin{equation}
\label{equationsGralY}
\frac{1}{2\pi}\, \mathfrak G\left( 2\pi \zeta_{k}^{(n)} \right) = \left\{ \frac{\theta_n}{2\pi}\right\}+k, \quad k\in \Z.
\end{equation}
In particular,
$$
[k,k+1) \ni \frac{1}{2\pi}\, \mathfrak G\left(  2\pi \zeta_{k}^{(n)} \right)=  \frac{1}{2\pi}\,\left(  2\pi \zeta_{k}^{(n)} - 2 \arg \left(  {_{1}F_{1}}\left( \lambda;1; 2\pi i \zeta_{k}^{(n)} \right)\right)\right)\geq \zeta_{k}^{(n)},
$$
where we have used \emph{(ii)} of Proposition \ref{prop:nozeros}. Hence,
$$
0\leq  \zeta_{0}^{(n)}<1 \quad \text{and} \quad  \zeta_{k-1}^{(n)}<  \zeta_{k}^{(n)}<k+1, \quad k\in \Z.
$$
By compactness and diagonal argument, we can always select a subsequence $\Lambda \subset \N$ such that the following limits exist:
$$
\lim_{n\in \Lambda} \zeta_{k}^{(n)}=\zeta_{k}, \quad k\in \Z.
$$
By \eqref{equationsGralY},
$$
\frac{1}{2\pi}\, \mathfrak G\left( 2\pi \zeta_{k}^{(n)} \right)-\frac{1}{2\pi}\, \mathfrak G\left( 2\pi \zeta_{k-1}^{(n)} \right)=1,
$$
and taking limits we conclude that
\begin{equation}
\label{eqforYbis}
\zeta_{k}-\zeta_{k-1}=1 + \frac{1}{\pi}\, \big( \arg \left(  {_{1}F_{1}}\left( \lambda;1; 2\pi i \zeta_{k} \right)\right)-  \arg \left(  {_{1}F_{1}}\left( \lambda;1; 2\pi i \zeta_{k-1} \right)\right)\big).
\end{equation}
Let $k\in \N$; since $ \arg \left(  {_{1}F_{1}}\left( \lambda;1; 2\pi i \zeta_{k} \right)\right)$ is strictly decreasing in $[0,+\infty)$,  the second term in the right hand side of \eqref{eqforYbis} is $< 0$, so that we conclude that
$$
0<\zeta_{k}-\zeta_{k-1}<1, \quad k\in \N.
$$
By \eqref{xandy}, we obtain that
$$
0<\liminf_{n} \frac{n }{\pi} \left( x^{(n)}_{k} -  x_{k-1}^{(n)} \right)
\leq \limsup_{n}\frac{n }{\pi} \left(x^{(n)}_{k} -  x_{k-1}^{(n)} \right)<1, \quad k \in \N.
$$

In the same vein, since $ \arg \left(  {_{1}F_{1}}\left( \lambda;1; 2\pi i \zeta_{k} \right)\right)$ is strictly increasing in $(-\infty,0)$,  by \eqref{eqforYbis},
$$
\zeta_{k}-\zeta_{k-1}>1, \quad -k\in \N,
$$
so that
$$
\liminf_{n} \frac{n }{\pi} \left( x^{(n)}_{k} -  x_{k-1}^{(n)} \right)>1, \quad -k \in \N.
$$

Furthermore, observe that for $c\neq 1$, the accumulation points of the sequence $\zeta_{0}^{(n)}$ is dense in the interval $ \mathfrak G^{-1}([0,2\pi])$. Indeed, by \eqref{defOfEta} and \eqref{defOfTheta},
$$
\frac{\theta_n}{2\pi}=\frac{n}{2}+ \frac{\log c}{\pi^{2}}\, \log n + \upsilon, \quad \upsilon \isdef \frac{\log c}{\pi^{2}}\, \log 4 +\frac{\beta-\alpha}{4}+ \frac{\hbar(0)}{\pi} -\frac{1}{\pi}\, \arg \Gamma(\lambda).
$$
Since $c\neq 1$, we can always take $b\in \{2, 3\}$ such that $(\log c)(\log b)/\pi^{2}\notin \mathbb Q$ (indeed, otherwise we would have that $\log 3/\log 2$ is rational, which is obviously impossible). For such a $b$, with $n=2\, b^{m}$,  $m\in \N$, equation \eqref{equationsGralY} is rewritten as
$$
\frac{1}{2\pi}\, \mathfrak G\left(2\pi  \zeta_{0}^{(n)} \right) = \left\{m \frac{\log c}{\pi^{2}}\, \log b + \upsilon + \frac{\log c}{\pi^{2}}\, \log 2 \right\}  .
$$
By Kronecker-Weyl theorem (see, e.g.~\cite[Chapter III]{Cassels:1957}),  the sequence
$$
\left\{m \frac{\log c}{\pi^{2}}\, \log 2 + \upsilon + \frac{\log c}{\pi^{2}}\, \log 2 \right\}
$$
is dense in $(0,1)$, and it remains to use the strict monotonicity of $\mathfrak G$. This finishes the proof of Proposition \ref{prop:clockbehavior}.

\begin{remark}
Obviously, if $c=1$, $\lambda=0$, ${_{1}F_{1}}\left(  \lambda;1;2\pi i y\right)\equiv 1$, and we obtain the clock behavior via \eqref{xandy} and \eqref{eqforYbis}.
\end{remark}

Regarding Remark \ref{rem:quasi}, it is easy to check numerically that in general
$$
\zeta_{k}-\zeta_{k-1}\neq \zeta_{k+1}-\zeta_{k}, \quad k\in \Z,
$$
which shows that the quasi-clock behavior fails too. Nevertheless, for function $y(x)$ introduced in Proposition \ref{prop:nozeros} we obtain applying \eqref{3} that $y'(x)=\mathcal O(1/|x|)$, $|x|\to\infty$, from where the ``clock behavior in the limit'' \eqref{clockinthelimit} follows by \eqref{eqforYbis}.

\subsection{Christoffel-Darboux Kernel}

Using the Christoffel-Darboux formula \cite[Section 3.2]{szego:1975}, we can write the kernel \eqref{defCDKernel} as
\begin{align*}
K_{n}\left(  x,y\right)   &  =\frac{k_{n-1}}{k_{n}}\frac{p_{n}\left(
x\right)  p_{n-1}\left(  y\right)  -p_{n}\left(  y\right)  p_{n-1}\left(
x\right)  }{x-y}\\
&  =k_{n-1}^{2}\frac{P_{n}\left(  x\right)  P_{n-1}\left(  y\right)
-P_{n}\left(  y\right) P_{n-1}\left(  x\right)  }{x-y}, \quad x\neq y;
\end{align*}
using  (\ref{sol Y}) and the fact that $\det \mathbf Y \equiv 1$, we obtain%
\begin{align}
K_{n}\left(  x,y\right)   &  =\frac{-1}{2\pi i}\frac{\mathbf Y_{11} \left(  x,n\right)  \mathbf Y_{21} \left(  y, n \right)  -\mathbf Y_{11}%
 \left(  y, n \right)  \mathbf Y_{21} \left(
x,n\right)  }{x-y}\nonumber \\
&  =\frac{1}{2\pi i}\frac{1}{x-y}
\begin{pmatrix}
0, & 1
\end{pmatrix} \mathbf  Y^{-1}\left(  y,n\right) \mathbf  Y\left(  x,n\right)  \begin{pmatrix}
1 \\ 0
\end{pmatrix} . \label{dem ass k 0}%
\end{align}

By analyticity, it is obviously sufficient to compute $K_n$ when $x, y\in (0,\delta )$. From \eqref{intermediateY},
\begin{equation} \label{Yinx}
  W(x) \mathbf Y_+\left(  x,n\right) \begin{pmatrix}
1\\
0
\end{pmatrix} =2^{-n\sigma_{3}}\mathbf R\left(  x\right)   D_{\infty}^{\sigma_{3}} \mathbf A_+(x)\,
 m_n(x)^{\sigma_3}\begin{pmatrix}
\Gamma\left(  1-\lambda\right)  G\left(  \lambda ;nf\left(  x\right)\right)   \\
\Gamma\left(  1+\lambda\right)  G\left(  1+\lambda ;nf\left(  x\right)\right)
\end{pmatrix}.
\end{equation}
On the other hand, since $\det \mathbf Y=1$, we have by \eqref{dem ass k 2}, \eqref{DefinitionEN}, \eqref{sol-N}, \eqref{def_m},
\begin{equation*}
\begin{split}
\mathbf Y\left(  y,n\right)^{-1}  &=  \varphi\left(  y\right)  ^{-n\sigma_{3}} \begin{pmatrix}
1 & 0\\
-\frac{1}{w_c}\, \varphi^{-2n} & 1
\end{pmatrix} \varphi\left(  y\right)  ^{n\sigma_{3}}W\left(  y\right)
^{\sigma_{3}}[ \mathbf \Psi\left(  nf\left(  y\right)  \right) ]^{-1} \\
&\times m_n(y)^{-\sigma_3}  \mathbf A_+(y)^{-1}
D_{\infty}^{-\sigma_{3}}
\mathbf R\left( y\right)   ^{-1}
2^{n\sigma_{3}} .
\end{split}
\end{equation*}
The matrix $\mathbf \Psi$ built in \eqref{sol Psi1}--\eqref{sol Psi6} also satisfies $\det\mathbf \Psi=1$, so that considerations that lead us to \eqref{intermediateY} show that
\begin{equation} \label{YinY}
\begin{split}
\begin{pmatrix}
0, & 1
\end{pmatrix}\, W(y) \mathbf Y\left(  y,n\right)^{-1}  &=  \begin{pmatrix}
- \Gamma\left(  1+\lambda\right)  G\left(  1+\lambda ;nf\left(  y\right)\right), \Gamma\left(  1-\lambda\right)  G\left(  \lambda ;nf\left(  y\right)\right)
\end{pmatrix} \\
&\times m_n(y)^{-\sigma_3}  \mathbf A_+(y)^{-1}
D_{\infty}^{-\sigma_{3}}
\mathbf R\left( y\right)   ^{-1}
2^{n\sigma_{3}} . 
\end{split}
\end{equation}
Observe also that locally uniformly for $z \in (-\delta, \delta)$, $\mathbf R(z)=\mathbf  I+ \OO(1/n)$, which implies that
$$
\mathbf R\left( y\right)   ^{-1} \mathbf R\left( x\right) = \mathbf  I+ \OO(1/n), \quad \text{locally uniformly for } x, y \in (-\delta, \delta).
$$
Gathering \eqref{Yinx} and \eqref{YinY} in \eqref{dem ass k 0} we conclude that for $x\neq y$,
\begin{align*}
\lim_{n\to \infty}\frac{\pi}{n}\, W\left( \frac{\pi x }{n}\right)  W\left( \frac{\pi y }{n}\right) K_{n}\left(  \frac{\pi x }{n},\frac{\pi y }{n}\right)   &    =\frac{1}{2\pi i}\frac{1}{ x-y } \\
\times
\begin{pmatrix}
- \Gamma\left(  1+\lambda\right)  G\left(  1+\lambda ;2 \pi i y \right), \Gamma\left(  1-\lambda\right)  G\left(  \lambda ;2 \pi i y \right)
\end{pmatrix} & \begin{pmatrix}
\Gamma\left(  1-\lambda\right)  G\left(  \lambda ;2 \pi i x\right)   \\
\Gamma\left(  1+\lambda\right)  G\left(  1+\lambda ;2 \pi i x \right)
\end{pmatrix}.
\end{align*}

By formula (6.1.31) of \cite{abramowitz/stegun:1972} and the definition of $\lambda$ in \eqref{def_lambda},
\begin{equation*}
 \Gamma\left(  1+ \lambda \right)\Gamma\left( 1-   \lambda \right)  =\frac{ \log c }{\sinh\left( \log c  \right)  } =\frac{2 c \log c}{c^2-1}. 
\end{equation*}
Since
$$
\lim_{n\to \infty}  W\left( \frac{\pi x }{n}\right)=\sqrt{c\, h(0)},
$$
we obtain \eqref{kernelFinal1}--\eqref{limitkernel}. Taking into account \eqref{conjugates} we can easily rewrite this formula in the form \eqref{kernelFinal2}.

Finally, the confluent form of the kernel in \eqref{limitkernel} is obtained from the expression for $x\neq y$ by taking limit $y\to x$.

\section{Properties of the confluent hypergeometric function} \label{sec:confluentproperties}

We prove finally the properties of ${_{1}F_{1}}$ and related functions summarized in  Proposition \ref{prop:nozeros}.

\emph{(i)} Let $K_{n}$ be the Christoffel-Darboux kernel defined in \eqref{defCDKernel}.  
Then for $z\in \C$,
$$
K_{n}(z,\overline{z})= \sum_{k=0}^{n-1}|p_{k}\left(  z\right)  |^{2} >0.
$$
This property is obviously inherited in the limit \eqref{kernelFinal1}, which implies that
$$
K_{\infty}(z,\overline{z}) \geq 0, \quad z\in \C.
$$
On the other hand, from formula (13.1.27) of \cite{abramowitz/stegun:1972} it follows that
\begin{align}
\label{conjugateF}
\overline{ {_{1}F_{1}}\left(  i a;1; i \overline{z}\right)} & = {_{1}F_{1}}\left(-i a;1; - i z\right)=e^{-z i} {_{1}F_{1}}\left(  ia+1;1;  i z\right), \\
\overline{G\left(  i a;  i \overline{z}\right)} & = G\left(  ia+1;   i z\right).
\label{conjugateG}
\end{align}
Thus, for $z\in \C\setminus \R$, $ \lambda=ia$, $a\in \R \setminus \{0\}$,
\begin{align} \nonumber
K_{\infty}\left( z,\overline{z} \right)  &= \dfrac{1}{2 \pi \, h(0) } \dfrac{  \log c}{c^2-1} \, \dfrac{    G\left(   \lambda ;2 \pi i z \right)  G\left( 1+  \lambda ;2 \pi i \overline{z} \right) - G\left(  1+\lambda ;2 \pi i z \right)  G\left(   \lambda ;2 \pi i \overline{z} \right) }{    \Im z } \\
& =  \dfrac{1}{2 \pi \, h(0) } \dfrac{  \log c}{c^2-1} \, \dfrac{   | G\left(   \lambda ;2 \pi i z \right)  |^{2} - |G\left(  \lambda ;2 \pi i \overline{z} \right) |^{2} }{    \Im z }\geq 0,
\label{inequalityForK}
\end{align}
which yields that 
\begin{equation}\label{inequality1forG}
    |G(\lambda ; i z)| \geq  |G(\lambda ; i \overline{z})|, \quad \Im z >0. 
\end{equation}
Assume that for $\zeta \in \C^+$, $G(\lambda ; i \zeta)=0$; then by \eqref{conjugateG} and \eqref{inequality1forG},
$$
G(\lambda ; i \overline{\zeta})=G\left(  ia+1;   i \zeta \right)=0.
$$
Hence, 
$$
{_{1}F_{1}}\left( i a;1;   i \zeta \right)= {_{1}F_{1}}\left( 1+i a;1;   i \zeta \right)=0.
$$
By induction and recurrence relation (13.4.1) in  \cite{abramowitz/stegun:1972} we conclude that every ${_{1}F_{1}}\left(  i a+n;1;  i \zeta \right)$, with $n\in \Z\cup \{ 0\}$, vanishes. But this is impossible, as follows from the addition formula
\begin{equation}\label{additionFormula}
    {_{1}F_{1}}\left( \lambda;1;  z+\zeta \right)=\left(\frac{\zeta}{z+\zeta} \right)^{\lambda} \sum_{n=0}^{\infty} \frac{(\lambda)_{n} \, z^{n}}{n!\, (z+\zeta )^{n}}\,  {_{1}F_{1}}\left(  \lambda+n;1;   \zeta\right)
\end{equation}
(see \cite[formula (2.3.4)]{Slater:1960je}).

Thus, we conclude that $f_{1}(z)=G(\lambda,   i z) \in \overline{HB}$. The assertion for $f_{2}$ is obtained by means of formula \eqref{conjugateF}.

\emph{(ii)}  
In order to prove \eqref{Fnonzero} assume that  $a\in \R\setminus \{0\}$ and $x\in \R$. Then by \eqref{conjugateF},
\begin{equation}
\label{conjugateF1}
\overline{ {_{1}F_{1}}\left(  i a;1; i x\right)}= {_{1}F_{1}}\left(-i a;1; - i x\right)=e^{-x i} {_{1}F_{1}}\left(  ia+1;1;  i x\right).
\end{equation}
Hence, an assumption that for $x\in \R$, ${_{1}F_{1}}\left(  i a;1;  i x\right)=0$ implies that ${_{1}F_{1}}\left(  i a+1;1;  i x\right)$ also vanishes, and we arrive at a contradiction reasoning as above and using the addition formula \eqref{additionFormula}.

Furthermore, the location of the zeros in the corresponding half planes and the inequality \eqref{inequalitystrictforF} is a direct consequence of \emph{(i)}.  In particular, function 
$$
h(z)= \frac{{_{1}F_{1}}\left( 1+ i a ;1;  i z\right)}{{_{1}F_{1}}\left(   i a ;1;  i z\right)}
$$
is holomorphic in $\C^{+}$, continuous in $\overline{\C^{+}}=\C^{+}\cup \R$, and satisfies $|h(z)|\leq 1$ for $z\in \overline{\C^{+}}$ and $|h(z)|=1$ for $z\in \R$. In consequence, by the maximum principle, $|h(z)|<1$ for $z\in  \C^{+}$, which proves that the inequality in \eqref{inequalitystrictforF} for $z$ in the upper half plane is strict.

\emph{(iii)} Due to \emph{(ii)}, $y(x)$ is correctly defined and real-analytic on $\R$, in particular, $y'$ can vanish only at a discrete set of points that can accumulate only at infinity. Again by  (13.1.27) of \cite{abramowitz/stegun:1972},
\begin{equation}
\label{identityAdri}
\frac{{_{1}F_{1}}\left(1+i a;1;  i x\right)}{{_{1}F_{1}}\left(i a;1; i x\right)}= e^{i x}\, \frac{{_{1}F_{1}}\left(-i a;1;  -i x\right)}{{_{1}F_{1}}\left(i a;1; i x\right)}=e^{i(x-2y(x))}\,.
\end{equation}
With the straightforward identity (see (13.4.4) in \cite{abramowitz/stegun:1972})
$$
{_{1}F_{1}}\left(1+i a;1;  i x\right)-{_{1}F_{1}}\left( i a;1;  i x\right)= (ix)\, {_{1}F_{1}}\left(1+i a;2;  i x\right)
$$
we rewrite \eqref{identityAdri} as
\begin{equation}
\label{decreasingArg1}
e^{i(x-2 y(x))} = 1 + ix\, \frac{{_{1}F_{1}}\left(1+i a;2;  i x\right)}{{_{1}F_{1}}\left(i a;1; i x\right)}.
\end{equation}
Since the real part of the left hand side is $\leq 1$, this implies that
\begin{equation}\label{decreasingArg2}
\Im \, \left( \frac{{_{1}F_{1}}\left(1+i a;2;  i x\right)}{{_{1}F_{1}}\left(i a;1; i x\right)}\right)\begin{cases} \geq 0, & \text{for } x>0, \\ \leq 0, & \text{for } x<0.\end{cases}
\end{equation}
On the other hand, by (13.4.8) of \cite{abramowitz/stegun:1972},
\begin{equation}
\label{4}
y'(x)=\Im \left( \frac{d}{dx} \log \bigg({_{1}F_{1}}\left(i a;1; i x\right)\bigg) \right)= -a \Im \, \left( \frac{{_{1}F_{1}}\left(1+i a;2;  i x\right)}{{_{1}F_{1}}\left(i a;1; i x\right)}\right),
\end{equation}
and by inequality \eqref{decreasingArg2}, the first part of the statement of \emph{(ii)} follows. The differential equation in \eqref{3} is obtained by taking the real part in \eqref{decreasingArg1} and using \eqref{4}.

In order to prove \emph{(iv)} we observe that $\mathfrak G$ satisfies the following initial value problem:
\begin{equation}
\label{ODEforG}
x \, \mathfrak G'(x)=x + 2a \left( 1-\cos \mathfrak G(x)\right),\quad \mathfrak G(0)=0.
\end{equation}
For $a=0$ the statement is trivial. Assume first that $a>0$; then by \emph{(ii)}, we only need to prove that $\mathfrak G'(x)>0$ for $x<0$.

Since $\mathfrak G$ is also real analytic, expanding it at $x=0$ we readily conclude from \eqref{ODEforG} that $\mathfrak G'(0)=1$. Hence, $\mathfrak G$ is locally increasing at the origin. Differentiating  \eqref{ODEforG} we obtain that
\begin{equation}
\label{ODEforG2}
x \, \mathfrak G''(x)=1 + \mathfrak G'(x) \left( -1+2a \sin \mathfrak G(x)\right).
\end{equation}
If for $x=\zeta<0$, $\mathfrak G'(\zeta)=0$, then by \eqref{ODEforG2},
\begin{equation}
\label{zeta1}
\zeta\, \mathfrak G''(\zeta)=1.
\end{equation}
In particular, $\mathfrak G''(\zeta)<0$, which shows that every critical point  of $\mathfrak G$ in the negative semi-axis is a strict local maximum, which is incompatible with the behavior at the origin. Thus, $\mathfrak G'$ is sign-invariant on $(-\infty, 0)$, and in consequence, $\mathfrak G'(x)>0$ there.

Assume now $a<0$; again by \emph{(ii)}, we only need to prove that $\mathfrak G'(x)>0$ for $x>0$.
Reasoning as above, if for $\zeta>0$ we have $\mathfrak G'(\zeta)=0$, then we get \eqref{zeta1}, which shows that every critical point  of $\mathfrak G$ in the positive semi-axis is a strict local minimum, which is again incompatible with the behavior at the origin.

\begin{remark}\label{remk:thanks}
The proof of \emph{(i)} presented here is an evolution of an idea of Doron Lubinsky to look at the Christoffel-Darboux kernels $K_{n}$ for each finite $n\in \N$ and then taking limits as $n\to \infty$; see also \cite{Lubinsky2009}.
The argument that allows to conclude that  ${_{1}F_{1}}$ does not vanish, based on the addition formula \eqref{additionFormula}, and the identity \eqref{identityAdri} were suggested to us by Adri Olde Daalhuis.
We gratefully acknowledge these two contributions of our colleagues.

\end{remark}

\section*{Acknowledgements}

AMF is partially supported by Junta de Andaluc\'{\i}a, grants FQM-229 and P06-FQM-01735, as well as a grant from the Ministry of Science and Innovation of Spain
(project code MTM2008-06689-C02-01).

VPS is sponsored by FCT (Portugal), under contract/grant
SFRH/BD/29731/2006.

Doron Lubinsky and Adri Olde Daalhuis contributed with elegant ideas to the proof of Proposition \ref{prop:nozeros} (see Remark \ref{remk:thanks} above).
We are grateful to Alexei Borodin for his interest and stimulating discussions, and in particular, for driving our attention to reference \cite{Borodin:2001xr}. Finally, the referees played an important role providing relevant suggestions and catching several typos, which  denotes their careful reading of such a technical manuscript.


\begin{thebibliography}{10}
\expandafter\ifx\csname url\endcsname\relax
  \def\url#1{\texttt{#1}}\fi
\expandafter\ifx\csname urlprefix\endcsname\relax\def\urlprefix{URL }\fi

\bibitem{abramowitz/stegun:1972}
M.~Abramowitz, I.~A. Stegun, Handbook of Mathematical Functions, Dover Publ.,
  New York, 1972.

\bibitem{MR2001m:05258a}
J.~Baik, P.~Deift, K.~Johansson, On the distribution of the length of the
  second row of a {Y}oung diagram under {P}lancherel measure, Geom. Funct.
  Anal. 10~(4) (2000) 702--731.

\bibitem{Borodin:2001xr}
A.~Borodin, G.~Olshanski, Infinite random matrices and ergodic measures,
Comm. Math. Phys. 223 (1) (2001) 87--123.

\bibitem{Cassels:1957}
J.~W.~S. Cassels, An introduction to {D}iophantine approximation, Cambridge
  Tracts in Mathematics and Mathematical Physics, No. 45, Cambridge University
  Press, New York, 1957.

\bibitem{MR2001f:42037}
P.~Deift, T.~Kriecherbauer, K.~T.-R. McLaughlin, S.~Venakides, X.~Zhou, Strong
  asymptotics of orthogonal polynomials with respect to exponential weights,
  Comm. Pure Appl. Math. 52~(12) (1999) 1491--1552.

\bibitem{MR98b:35155}
P.~Deift, S.~Venakides, X.~Zhou, New results in small dispersion {K}d{V} by an
  extension of the steepest descent method for {R}iemann-{H}ilbert problems,
  Internat. Math. Res. Notices~(6) (1997) 286--299.

\bibitem{MR94d:35143}
P.~Deift, X.~Zhou, A steepest descent method for oscillatory
  {R}iemann-{H}ilbert problems. {A}symptotics for the {M}{K}d{V} equation, Ann.
  of Math. 137~(2) (1993) 295--368.

\bibitem{MR2000g:47048}
P.~A. Deift, Orthogonal polynomials and random matrices: a {R}iemann-{H}ilbert
  approach, New York University Courant Institute of Mathematical Sciences, New
  York, 1999.

\bibitem{MR96d:34004}
P.~A. Deift, X.~Zhou, Asymptotics for the {P}ainlev\'e {I}{I} equation, Comm.
  Pure Appl. Math. 48~(3) (1995) 277--337.

\bibitem{Fokas92}
A.~Fokas, A.~Its, A.~Kitaev, The isomonodromy approach to matrix models in {2D}
  quantum gravity, Comm. Math. Phys. 147 (1992) 395--430.

\bibitem{FMS2}
A.~Foulqui\'e Moreno, A.~Mart\'{\i}nez-Finkelshtein, V.L.~Sousa, On a Conjecture of A. Magnus concerning the asymptotic behavior of the recurrence coefficients of the generalized Jacobi polynomials, Journal of Approximation Theory, in press. Also preprint arXiv:0905.2753.

\bibitem{Gakhov:90}
F.~D. Gakhov, Boundary value problems, Dover Publications Inc., New York, 1990,
  translated from the Russian, Reprint of the 1966 translation.

\bibitem{Its07b}
A.~Its, I.~Krasovsky, Hankel determinant and orthogonal polynomials for the
  gaussian weight with a jump, Contemp. Math. 458 (2008) 215--247.

\bibitem{MR2087231}
A.~B.~J. Kuijlaars, K.~T.-R. McLaughlin, W.~Van~Assche, M.~Vanlessen, The
  {R}iemann-{H}ilbert approach to strong asymptotics for orthogonal polynomials
  on {$[-1,1]$}, Adv. Math. 188~(2) (2004) 337--398.

\bibitem{Lubinsky2009}
D.~S. Lubinsky, Universality limits for random matrices and de Branges spaces of entire functions, Journal of Functional Analysis 256 (2009) 3688--3729.

\bibitem{Lubinsky2008}
E.~Levin, D.~S. Lubinsky, Applications of universality limits to zeros and
  reproducing kernels of orthogonal polynomials, Journal of Approximation
  Theory 150 (2008) 69--95.

\bibitem{Magnus1995}
A.~P. Magnus, Asymptotics for the simplest generalized {J}acobi polynomials
  recurrence coefficients from {F}reud's equations: numerical explorations,
  Ann. Numer. Math. 2 (1995) 311--325.

\bibitem{Simon2008}
B.~Simon, The {C}hristoffel-{D}arboux kernel,
in ``Perspectives in PDE, Harmonic Analysis and Applications,'' a volume in honor of V.G. Maz'ya's 70th birthday, Proceedings of Symposia in Pure Mathematics 79 (2008) 295--335.

\bibitem{Simon:2009ly}
B.~Simon, Fine structure of the zeros of orthogonal polynomials: a progress
report, preprint (2009), to appear in ``Recent Trends in Orthogonal
Polynomials and Approximation Theory'',  a volume in honor of G.~L\'opez Lagomasino's 60th birthday,
Contemporary Mathematics series, AMS.

\bibitem{Slater:1960je}
L.~J. Slater, Confluent hypergeometric functions, Cambridge University Press,
  Cambridge, UK, 1960.

\bibitem{szego:1975}
G.~Szeg{\H{o}}, Orthogonal Polynomials, vol.~23 of Amer.\ Math.\ Soc.\ Colloq.\
  Publ., 4th ed., Amer.\ Math.\ Soc., Providence, {RI}, 1975.

\end{thebibliography}
%
\def\cprime{$'$}

\medskip

\obeylines
\texttt{A.~Foulqui\'{e} Moreno (foulquie@ua.pt)
Department of Mathematics and Research Unity "Matem\'{a}tica e Aplica\c{c}\~{o}es",
University of Aveiro, Campus Universit\'ario de Santiago,
3810-193 Aveiro, Portugal
\medskip
A. Mart\'{\i}nez-Finkelshtein (andrei@ual.es)
Department of Statistics and Applied Mathematics
University of Almer\'{\i}a, Spain, and
Instituto Carlos I de F\'{\i}sica Te\'{o}rica y Computacional
Granada University, Spain
\medskip
V.L. Sousa (vsousa@ua.pt)
Escola Secund\'{a}ria Jo\~{a}o Silva Correia, and 
Research Unit "Matem\'{a}tica e Aplica\c{c}\~{o}es",
University of Aveiro, Campus Universit\'ario de Santiago,
3810-193 Aveiro, Portugal
}

\end{document}